\documentclass[12pt]{amsart}

\newtheorem{cotmb}{Corollary}
\newtheorem{prop}{Proposition}
\newcommand{\cal}{\mathcal}

\newtheorem{lemma}{Lemma}

\usepackage[mathscr]{euscript}

\vspace{0.5cm}

\title[Harmonicity of almost complex structures]
{Harmonicity of the Atiyah-Hitchin-Singer and Eells-Salamon almost complex structures}
\author{Johann Davidov, Oleg Mushkarov}

\thanks{The authors are partially supported by  the National Science
Fund, Ministry of Education and Science of Bulgaria under contract
DFNI-I 02/14. }

\address{Johann Davidov\\Institute of Mathematics and Informatics \\
Bulgarian Academy of Sciences\\ Acad. G.Bonchev Str. Bl.8\\ 1113
Sofia\\ Bulgaria\\  \newline \centerline{and} \newline University of
Structural Engineering and Architecture "L.Karavelov", 175
Suhodolska St., 1373 Sofia, Bulgaria} \email{jtd@math.bas.bg}

\address{Oleg Mushkarov \\Institute of Mathematics and Informatics \\
Bulgarian Academy of Sciences\\ Acad. G.Bonchev Str. Bl.8\\ 1113
Sofia\\ Bulgaria\\ \newline \centerline{and} \newline Sauth-West
University, 2700 Blagoevgrad, Bulgaria}
\email{muskarov@math.bas.bg}

\begin{document}

\begin{abstract}
In this paper we describe the oriented Riemannian four-manifolds $M$
for which the Atiyah-Hitchin-Singer or Eells-Salamon almost complex
structure on the twistor space ${\mathcal Z}$ of $M$ determines a
harmonic map from ${\mathcal Z}$ into its twistor space.

\vspace{0,1cm} \noindent 2010 {\it Mathematics Subject
Classification}. Primary 53C43, Secondary 58E20, 53C28

\vspace{0,1cm} \noindent {\it Key words: twistor spaces, almost
complex structures, harmonic maps}
\end{abstract}

\thispagestyle{empty}

\maketitle


\section{Introduction}\label{Intro}

The twistor approach has been used for years for studying
conformal geometry of four-manifolds  by means of complex
geometric methods, and in this way many important results have
been obtained. Moreover, the twistor spaces endowed with the
Atiyah-Hitchin-Singer and Eells-Salamon almost complex structures
are interesting geometric objects in their own right whose
geometric properties have been studied by many authors. In this
paper we look at these structures from the point of view of
variational theory. The motivation behind is the fact that if a
Riemannian manifold admits an almost complex structure compatible
with its metric, it possesses many such structures (cf., for
example, \cite{D16, DHM15}). Thus,  it is natural to seek criteria
that distinguish some of these structures among all. One way to
obtain such a criterion is to consider the compatible almost
complex structures on a Riemannian manifold $(N,h)$ as sections of
its twistor bundle ${\mathcal Z}$. The smooth manifold ${\mathcal
Z}$ admits a natural Riemannian metric $h_1$ such that the
projection map $({\mathcal Z},h_1)\to (N,h)$ is a Riemannian
submersion with totally geodesic fibres.  From this point of view,
E. Calabi and H. Gluck \cite{CG} have proposed  to consider as
"the best" those compatible almost complex structures $J$ on
$(N,h)$ whose image $J(N)$ in ${\mathcal Z}$ is of minimal volume.
They have proved that the standard almost Hermitian structure on
the $6$-sphere $S^6$, defined by means of the Cayley numbers, can
be characterized by that property. Another criterion has been
discussed by C. M. Wood \cite{W1,W2} who has suggested to single out
the structures $J$ that are harmonic sections of the twistor
bundle ${\mathcal Z}$, i.e. critical points of the energy
functional under variations through sections of ${\mathcal Z}$.
While the K\"ahler structures are absolute minima of the energy
functional, there are many examples of non-K\"ahler structures,
which are harmonic sections \cite{W1,W2}. Sufficient conditions
for a compatible almost complex structure to be a minimizer of the
energy functional and examples of non-K\"ahler minimizers have
been given by G. Bor, L. Hern\'andez-Lamoneda and M. Salvai
\cite{BLS}.

Forgetting the bundle structure of ${\mathcal Z}$, we can also
consider compatible almost complex structures that are critical
points of the energy functional under variations through all maps
$N\to{\mathcal  Z}$. These structures are genuine harmonic maps
from $(N,h)$ into $({\mathcal Z},h_1)$; we refer to \cite{EL} for
basic facts about harmonic maps. The problem when a compatible
almost complex structure on a four-dimensional Riemannian manifold
is a harmonic map into its twistor space has been studied in
\cite{DHM15} (see also \cite{D16}).

If the base manifold $N$ is oriented, the twistor space ${\mathcal Z}$
has two connected components often called positive and negative
twistor spaces of $(N,h)$; their sections are compatible  almost
complex structures yielding the orientation and, respectively, the
opposite orientation of $N$.

Setting $h_t=\pi^{\ast}h+th^v$,  $t>0$,  where $\pi:{\cal Z}\to N$ is the projection
map and $h^v$ is the metric of the fibre,  define a $1$-parameter family  of Riemannian
metrics on ${\cal Z}$ compatible with the almost complex structures $\cal J_1$ and
$\cal J_2$ on ${\mathcal Z}$ introduced respectively
by Atiyah-Hitchin-Singer \cite{AHS} and Eells-Salamon \cite{ES}.
In \cite{DM02} we have found geometric
conditions on an oriented four-dimensional Riemannian manifold under
which the almost complex structures $\cal J_1$ and $\cal J_2$ on its
negative twistor space $({\mathcal Z},h_t)$ are harmonic sections.

\smallskip

\noindent {\bf Theorem 1.}~ {\it Let $(M,g)$ be an oriented
Riemannian $4$-manifold and let $({\cal Z},h_t)$ be its negative
twistor space. Then:

\noindent $(i)$ The Atiyah-Hitchin-Singer almost complex structure
${\cal J}_1$ on $({\cal Z},h_t)$ is a harmonic section if and only
if $(M,g)$ is a self-dual manifold.

\noindent $(ii)$ The Eells-Salamon almost complex structure ${\cal
J}_2$ on $({\cal Z},h_t)$ is a harmonic section if and only if
$(M,g)$ is a self-dual manifold with constant scalar curvature}.

\smallskip

By a theorem of Atiyah-Hitchin-Singer \cite{AHS} the self-duality of
$(M,g)$ is a necessary and sufficient condition for the
integrability of the almost complex structure ${\cal J}_1$. In
contrast, the almost complex structure ${\cal J}_2$ is never
integrable by a result of Eells-Salamon \cite{ES} but it is very
useful for constructing harmonic maps.

\smallskip

The aim of the present paper is to find the four-manifolds for which
the almost complex structures $\cal J_1$ and $\cal J_2$ are harmonic
maps. More precisely, we prove the following

\smallskip

\noindent {\bf Theorem 2.}~ {\it Let ${\cal J}_1$ and ${\cal J}_2$ be the
Atiyah-Hitchin-Singer and  Eells-Salamon almost complex
structures on the (negative) twistor space $({\cal Z},h_t)$ of an
oriented Riemannian  four-manifold $(M,g)$. Each ${\cal J}_k$ (k=1 or 2) is a harmonic map if
and only if $(M,g)$ is either self-dual and Einstein, or is
locally the product of an open interval in ${\mathbb R}$ and a
$3$-dimensional Riemannian manifold of constant curvature}.

\smallskip

Note that any compact  self-dual Einstein manifold with positive
scalar curvature is isometric to the $4$-sphere $\mathbb S^4$ or  the
complex projective space $\mathbb {CP}^2$ with their standard
metrics \cite{FK, H81} (see also \cite[Theorem 13.30]{Besse}). In the case of negative scalar curvature, a
complete classification is not available yet and the only known
compact examples are quotients of the unit ball in $\mathbb C^2$
with the metric of constant negative curvature or the Bergman
metric. In contrast, there are many local examples of self-dual
Einstein metrics with a prescribed sign of the scalar curvature
(cf., e.g., \cite {Der, H95, Le1, Le2, Le3, Ped, Tod}). Note also
that every Riemannian manifold that locally  is the product of an
open interval in ${\mathbb R}$ and a $3$-dimensional Riemannian
manifold of constant curvature $c$ is locally conformally flat with
constant scalar curvature $6c$. It is not Einstein unless $c=0$,
i.e. Ricci flat.

The proof of Theorem 2  is based on an explicit formula for the
second fundamental form $\widetilde\nabla J_{\ast}$ of a compatible
almost complex structure $J$ on a Riemannian manifold considered as
a map from the manifold into its twistor space (Proposition 1). In
particular, it follows from Theorem 1 mentioned above that if the
vertical part of $\mathit{Trace}\widetilde\nabla {\cal J}_{k\,\ast}$
vanishes then the manifold $(M, J)$ is self-dual. This simplifies
the formulas for the values of the horizonal part of
$\mathit{Trace}\widetilde\nabla {\cal J}_{k\,\ast}$ at vertical and
horizontal vectors (Lemmas 1 and 2). Using these formulas we show
that the Ricci tensor of $(M,g)$ is parallel and three of its
eigenvalues coincide. Thus either $(M,g)$ is Einstein or exactly
three of the eigenvalues coincide. In the second case a result in
\cite[Lemma 1]{DGM} (essentially due to LeBrun and Apostolov)
implies that the simple eigenvalue vanishes, thus $(M,g)$  is
locally the product of an interval in $\mathbb R$ and a $3$-manifold
of constant curvature.

Note also that if $(h_t,{\cal J}_1)$ is a K\"ahler structure, then ${\cal J}_1$ is a totally geodesic map.
It is a result of Friedrich-Kurke \cite{FK} that $(h_t,{\cal J}_1)$ is K\"ahler exactly when the base manifold is self-dual and Einstein with positive scalar curvature $12/t$.  The necessary and sufficient conditions for ${\cal J}_1$ and ${\cal J}_2$ to be totally geodesic maps will be discussed elsewhere.

\smallskip

\noindent {\bf Acknowledgment}. We would like to thank the referee  whose remarks helped to improve the final version of the paper.

\section{Preliminaries}\label{prel}

\subsection{The manifold of compatible linear complex
structures}\label{CLCS}

Let $V$ be a real vector space of even dimension $n=2m$ endowed
with an Euclidean metric $g$. Denote by $F(V)$ the set of all
complex structures on $V$ compatible with the metric $g$, i.e.
$g$-orthogonal. This set has the structure of an imbedded
submanifold of the vector space $\mathfrak{so}(V)$ of skew-symmetric
endomorphisms of $(V,g)$.

The group $O(V)$ of orthogonal transformations of $(V,g)$ acts
smoothly and transitively on the set $F(V)$ by conjugation. The
isotropy subgroup at a fixed $J\in F(V)$ consists of the
orthogonal transformations commuting with $J$. Therefore $F(V)$
can be identified with the homogeneous space $O(2m)/U(m)$. In
particular, $dim\,F(V)=m^2-m$. Moreover, $F(V)$ has two connected
components. If we fix an orientation on $V$, these components
consist of all complex structures on $V$ compatible with the
metric $g$ and inducing $\pm$ the orientation of $V$; each of them
has the homogeneous representation $SO(2m)/U(m)$.

The tangent space of $F(V)$ at a point $J$ consists of all
endomorphisms $Q\in \mathfrak{so}(V)$ anti-commuting with $J$ and we have the
decomposition
\begin{equation}\label{decom}
\mathfrak{so}(V)=T_JF(V)\oplus \{S\in\mathfrak{so}(V):~ SJ-JS=0\}.
\end{equation}
This decomposition is orthogonal with respect to the metric $G(A,B)=-\frac{1}{n}\mathit{Trace} AB$ of $\mathfrak{so}(V)$
(the factor $1/n$ is chosen so that every $J\in F(V)$ to have unit norm). The metric $G$ on $F(V)$ is compatible with the
almost complex structure ${\cal J}$ defined by
$$
{\cal J}Q=JQ~\mbox{for}~ Q\in T_JF(V).
$$

Let $J\in F(V)$ and let $e_1,...,e_{2m}$ be an orthonormal basis
of $V$ such that $Je_{2k-1}=e_{2k}$, $k=1,...,m$. Define
skew-symmetric endomorphisms $S_{a,b}$, $a,b=1,...,2m$, of $V$
setting
$$
S_{a,b}e_c=\displaystyle{\sqrt{\frac{n}{2}}}(\delta_{ac}e_b-\delta_{bc}e_a),\quad
c=1,...,2m.
$$
The maps $S_{a,b}$, $1\leq a<b\leq 2m$, constitute a $G$-orthonormal
basis of $\mathfrak{so}(V)$.  Set
$$
\begin{array}{c}
A_{r,s}=\frac{1}{\sqrt 2}(S_{2r-1,2s-1}-S_{2r,2s}),\quad
B_{r,s}=\frac{1}{\sqrt 2}(S_{2r-1,2s}+S_{2r,2s-1}),\\[6pt]
r=1,...,m-1,\>s=r+1,...,m.
\end{array}
$$
Then $\{A_{r,s},B_{r,s}\}$ is a $G$-orthonormal basis of $T_JF(V)$
with $B_{r,s}={\cal J}A_{r,s}$.

Denote by $D$  the Levi-Civita connection of the metric $G$ on
$F(V)$. Let $X,Y$ be vector fields on $F(V)$ considered as
$\mathfrak{so}(V)$-valued functions on $\mathfrak{so}(V)$. By the Koszul formula, for
every $J\in F(V)$,
\begin{equation}\label{coder}
(D_{X}Y)_J=\frac{1}{2}(Y'(J)(X_J)+J\circ Y'(J)(X_J)\circ J)
\end{equation}
where $Y'(J)\in Hom(\mathfrak{so}(V),\mathfrak{so}(V))$ is the derivative of the function
$Y:\mathfrak{so}(V)\to \mathfrak{so}(V)$ at the point $J$. The latter formula easily
implies that $(G,{\cal J})$ is a K\"ahler structure on $F(V)$. Note
also that the metric $G$ is Einstein with scalar curvature
$\frac{m}{2}(m-1)(m^2-m)$ (see, for example,\cite{D05}).

\subsection {The four-dimensional case} Suppose that $dim\,V=4$.
Then, as is well-known, each of the two connected components of
$F(V)$ can be identified with the unit sphere $S^2$. It is often
convenient to describe this identification in terms of the space
$\Lambda^2V$. The metric $g$ of $V$ induces a metric on $\Lambda^2V$
given by
$$
g(x_1\wedge x_2,x_3\wedge
x_4)=\frac{1}{2}[g(x_1,x_3)g(x_2,x_4)-g(x_1,x_4)g(x_2,x_3)],
$$
the  factor $1/2$ being chosen in consistence with \cite{DM91,DM02}.
Consider the isomorphisms $\mathfrak{so}(V)\cong \Lambda^2V$ sending
$\varphi\in \mathfrak{so}(V)$ to the $2$-vector $\varphi^{\wedge}$ for which
$$
2g(\varphi^{\wedge},x\wedge y)=g(\varphi x,y),\quad x,y\in V.
$$
This isomorphism is an isometry with respect to the metric $G$ on
$\mathfrak{so}(V)$ and the metric $g$ on $\Lambda^2V$. Given $a\in
\Lambda^2V$, the skew-symmetric endomorphism of $V$ corresponding
to $a$ under the inverse isomorphism will be denoted by $K_a$.

Fix an orientation on $V$ and denote by $F_{\pm}(V)$ the set of
complex structures on $V$ compatible with the metric $g$ and
inducing $\pm$ the orientation of $V$. The Hodge star operator
defines an endomorphism $\ast$ of $\Lambda^2V$ with $\ast^2=Id$.
Hence we have the decomposition
$$
\Lambda^2V=\Lambda^2_{-}V\oplus\Lambda^2_{+}V
$$
where $\Lambda^2_{\pm}V$ are the subspaces of $\Lambda^2V$
corresponding to the $(\pm 1)$-eigen-values of the operator
$\ast$. Let $(e_1,e_2,e_3,e_4)$ be an oriented orthonormal basis
of $V$. Set
\begin{equation}\label{s-basis}
s_1^{\pm}=e_1\wedge e_2 \pm e_3\wedge e_4, \quad s_2^{\pm}=e_1\wedge
e_3\pm e_4\wedge e_2, \quad s_3^{\pm}=e_1\wedge e_4\pm e_2\wedge
e_3.
\end{equation}
Then $(s_1^{\pm},s_2^{\pm},s_3^{\pm})$ is an orthonormal basis of
$\Lambda^2_{\pm}V$. Note that this basis defines an orientation on
$\Lambda^2_{\pm}V$, which does not depend on the choice of the basis
$(e_1,e_2,e_3,e_4)$  (see, for example, \cite{D16}). We call this
orientation "canonical".

It is easy to see that the isomorphism $\varphi\to \varphi^{\wedge}$
identifies $F_{\pm}(V)$ with the unit sphere $S(\Lambda^2_{\pm}V)$
of the Euclidean vector space $(\Lambda^2_{\pm}V,g)$. Under this
isomorphism, if $J\in F_{\pm}(V)$, the tangent space
$T_JF(V)=T_JF_{\pm}(V)$ is identified with the orthogonal complement
$({\Bbb R} J^{\wedge})^{\perp}$ of the space ${\Bbb R}J^{\wedge}$ in
$\Lambda^2_{\pm}V$.

Consider the $3$-dimensional Euclidean space $(\Lambda^2_{\pm}V,g)$
with its canonical orientation and denote by $\times$ the usual
vector-cross product in it. Then if $a,b\in\Lambda^2_{\pm}V$, the
isomorphism $\Lambda^2V\cong \mathfrak{so}(V)$ sends $a\times b$ to
$\pm\frac{1}{2}[K_a,K_b]$. Thus, if $J\in F_{\pm}(V)$ and $Q\in
T_JF(V)=T_JF_{\pm}(V)$, we have
\begin{equation}\label{calJ}
({\cal J}Q)^{\wedge}=\pm (J^{\wedge}\times Q^{\wedge}).
\end{equation}

\subsection{The twistor space of an even-dimensional Riemannian
manifold}

Let $(N,g)$ be a  Riemannian manifold of dimension $n=2m$. Denote
by $\pi:{\cal Z}\to N$ the bundle over $N$ whose fibre at every
point $p\in N$ consists of all compatible complex structures on
the Euclidean vector space $(T_pN,g_p)$. This is the associated
bundle
$$
{\cal Z}=O(N)\times_{O(n)}F({\Bbb R^n})
$$
where $O(N)$ is the principal bundle of orthonormal frames on $N$
and $F({\Bbb R^n})$ is the manifold of complex structures on ${\Bbb
R^n}$ compatible with its standard metric. The manifold ${\cal Z}$
is called the twistor space of $(N,g)$.

The Levi-Civita connection of $(N,g)$ gives rise to a splitting
${\cal V}\oplus {\cal H}$ of the tangent bundle of any bundle
associated to $O(N)$ into vertical and horizontal parts. This
allows one to define a natural $1$-parameter family of Riemannain
metrics $h_t$, $t>0$, on the manifold ${\cal Z}$ sometimes called
"the canonical variation of the metric of $N$" (\cite[Chapter 9
G]{Besse}). For every $J\in {\cal Z}$, the horizontal subspace
${\cal H}_J$ of $T_J{\cal Z}$ is isomorphic via the differential
$\pi_{\ast J}$ to the tangent space $T_{\pi(J)}N$ and the metric
$h_t$ on ${\cal H}_J$ is the lift of the metric $g$ on
$T_{\pi(J)}N$, $h_t|{\cal H}_J=\pi^{\ast} g$. The vertical
subspace ${\cal V}_J$ of $T_J{\cal Z}$ is the tangent space at $J$
to the fibre of the bundle ${\cal Z}$ through $J$ and $h_t|{\cal
V}_J$ is defined as $t$ times the metric $G$ of this fibre.
Finally, the horizontal space ${\cal H}_J$ and the vertical space
${\cal V}_J$ are declared to be orthogonal. Then, by  the Vilms
theorem \cite{V}, the  projection $\pi:({\cal Z},h_t)\to (N,g)$ is
a Riemannian submersion with totally geodesic fibres (this can
also be proved directly).

The manifold ${\cal Z}$ admits two almost complex structures
${\cal J}_1$ and ${\cal J}_2$  defined in the case $dim\,N=4$ by
Atiyah-Hitchin-Singer \cite{AHS} and Eells-Salamon \cite{ES},
respectively.  On a vertical space ${\cal V}_J$, ${\cal J}_1$ is
defined to be the complex structure ${\cal J}_J$ of the fibre
through $J$, while ${\cal J}_2$ is defined as the conjugate
complex structure, i.e. ${\cal J}_{2}|{\cal V}_J=-{\cal J}_J$. On
a horizontal space ${\cal H}_J$, ${\cal J}_1$ and ${\cal J}_2$ are
both defined to be the lift to ${\cal H}_J$ of the endomorphism
$J$ of $T_{\pi(J)}N$. The almost complex structures ${\cal J}_1$
and ${\cal J}_2$ are compatible with each metric $h_t$.

\smallskip

Consider ${\cal Z}$ as a submanifold of the bundle
$$\pi:A(TN)=O(N)\times_{O(n)}so(n)\to N$$ of skew-symmetric
endomorphisms of $TN$. The inclusion of ${\cal Z}$ into $A(TN)$ is
fibre-preserving and, for every $J\in{\cal Z}$,  the horizontal
subspace ${\cal H}_J$ of $T_J{\cal Z}$ coincides with the horizontal
subspace of $T_JA(TN)$ since the inclusion of $F({\Bbb R^n})$ into
$so(n)$ is $O(n)$-equivariant.

The Levi-Civita connection of $(N,g)$ determines a connection on the
bundle $A(TN)$, both denoted by $\nabla$, and the corresponding
curvatures are related by
$$
(R(X,Y)\varphi)(Z)=R(X,Y)\varphi(Z)-\varphi(R(X,Y)Z)
$$
for $\varphi\in A(TN)$, $X,Y,Z\in TN$. The  curvature operator
${\cal R}$ is the self-adjoint endomorphism of $\Lambda ^{2}TN$
defined by
$$
 g({\cal R}(X\land Y),Z\land T)=g(R(X,Y)Z,T),\quad X,Y,Z,T\in TN.
$$
Let us note that we adopt the following definition for the curvature
tensor $R$ : $R(X,Y)=\nabla_{[X,Y]}-[\nabla_{X},\nabla_{Y}]$.

\smallskip

Let $(U,x_1,...,x_n)$ be a local coordinate system of $N$ and
$E_1,...,E_n$ an orthonormal frame of $TN$ on $U$. Define sections
$S_{ij}, 1\leq i,j\leq n$, of $A(TN)$ by the formula
\begin{equation}\label{Sij}
S_{ij}E_l=\sqrt{\frac{n}{2}}(\delta_{il}E_j - \delta_{lj}E_i),\quad
l=1,...,n.
\end{equation}
 Then $S_{ij}, i<j,$ form an orthonormal frame of $A(TN)$
with respect to the metric
$G(a,b)=\displaystyle{-\frac{1}{n}}\mathit{Trace}(a\circ b)\, ;a,b\in
A(TN)$. Set
$$\tilde x_{i}(a)=x_{i}\circ\pi(a),\, y_{jl}(a)=\sqrt{\frac{2}{n}}G(a,S_{jl}), j<l,$$ for
$a\in A(TN)$. Then $(\tilde x_{i},y_{jl})$ is a local coordinate
system of the manifold $A(TN)$.  Setting $y_{lk}=-y_{kl}$ for $l\geq
k$, we have $aE_j=\sum_{l=1}^n y_{jl}E_l$, j=1,...,n.

For each vector field
$$X=\textstyle{\sum\limits_{i=1}^n} X^{i}\frac{\partial}{\partial x_i}$$
on $U$, the horizontal lift $X^h$ on $\pi^{-1}(U)$ is given by
\begin{equation}\label{eq 3.1}
X^{h}=\textstyle{\sum\limits_{i=1}^n
}(X^{i}\circ\pi)\frac{\partial}{\partial\tilde x_i}-
\textstyle{\sum\limits_{j<l}\sum\limits_{p<q}}
y_{pq}G(\nabla_{X}S_{pq},S_{jl})\circ\pi\frac{\partial}{\partial
y_{jl}}.
\end{equation}

     Let $a\in A(TN)$ and $p=\pi(a)$. Then (\ref{eq 3.1}) implies
that, under the standard identification $T_{a}A(TN)\cong A(T_{p}N)$
(=the skew-symmetric endomorphisms of $(T_{p}N,g_p)$), we have
\begin{equation}\label{Lie-2}
[X^{h},Y^{h}]_{a}=[X,Y]^h_a + R(X,Y)a.
\end{equation}

\smallskip

 Farther we shall often make use of the isomorphism $A(TN)\cong\Lambda^{2}TN$ that assigns to each $a\in
A(T_{p}N)$ the 2-vector $a^{\wedge}$ for which
$$
2g(a^{\wedge},X\wedge Y)=g(aX,Y),\quad X,Y\in T_{p}N,
$$
the metric on $\Lambda^{2}TN$ being defined by
$$
g(X_1\wedge X_2,X_3\wedge X_4)=
\frac{1}{2}[g(X_1,X_3)g(X_2,X_4)-g(X_1,X_4)g(X_2,X_3)].
$$

\begin{lemma}\label{R[a,b]} {\rm (\cite{D05})} For every $a,b\in A(T_pN)$ and $X,Y\in T_pN$, we
have
\begin{equation}\label{Rw}
G(R(X,Y)a,b)=\frac{2}{n}g(R([a,b]^{\wedge})X,Y).
\end{equation}
\end{lemma}
{\bf Proof}. Let $E_1,...,E_n$ be an orthonormal basis of $T_pN$.
Then
$$
[a,b]=\frac{1}{2}\textstyle{\sum\limits_{i,j=1}^n} g([a,b]E_i,E_j)E_i\wedge E_j.
$$
Therefore
$$
\begin{array}{l}
g(R([a,b]^{\wedge})X,Y)\\[6pt]
=\frac{1}{2}\sum\limits_{i,j=1}^n
g(R(X,Y)E_i,E_j)[g(abE_i,E_j)+g(aE_i,bE_j)]\\[6pt]
=\frac{1}{2}\sum\limits_{i=1}^ng(R(X,Y)E_i,abE_i)\\[6pt]
\hspace{2.5cm}+\frac{1}{2}\sum\limits_{i,j,k=1}^n
g(R(X,Y)E_i,E_j)g(E_i,aE_k)g(E_j,bE_k)\\[6pt]
=-\frac{1}{2}\sum\limits_{i=1}^n
g(a(R(X,Y)E_i),bE_i)+\frac{1}{2}\sum\limits_{k=1}^n
g(R(X,Y)aE_k,bE_k)\\[6pt]
=\frac{n}{2}G(R(X,Y)a,b).
\end{array}
$$

For every $J\in{\cal Z}$, we identify the vertical space ${\cal
V}_J$ with the subspace of $A(T_{\pi(J)}N)$ of skew-symmetric
endomorphisms anti-commuting with $J$. Then, for every section $K$
of the twistor space ${\cal Z}$ near a point $p\in N$ and every
$X\in T_pN$, the endomorphism $\nabla_{X}K$ of $T_pN$ belongs to the
vertical space ${\cal V}_{K(p)}$.


Lemma~\ref{R[a,b]} implies that
\begin{equation}\label{r-r}
h_t(R(X,Y)J,V)=\frac{2t}{n}g(R([J,V]^{\wedge})X,Y)=
\frac{4t}{n}g(R((J\circ V)^{\wedge})X,Y).
\end{equation}

\smallskip

Denote by $D$ the Levi-Civita connection of $({\cal Z},h_t)$.

\begin{lemma}\label{LC} {\rm (\cite{D05,DM91})} If $X,Y$ are vector
fields on $N$ and $V$ is a vertical vector field on ${\cal Z}$, then
\begin{equation}\label{D-hh}
(D_{X^h}Y^h)_{J}=(\nabla_{X}Y)^h_{J}+\frac{1}{2}R_{p}(X\wedge Y)J
\end{equation}
\begin{equation}\label{D-vh}
(D_{V}X^h)_{J}={\cal H}(D_{X^h}V)_{J}=-\frac{2t}{n}(R_{p}((J\circ
V_J)^{\wedge})X)_J^h
\end{equation}
where $J\in{\cal Z}$, $p=\pi(J)$, and ${\cal H}$ means "the
horizontal component".
\end{lemma}
{\bf Proof}. Identity (\ref{D-hh}) follows from the Koszul formula
for the Levi-Civita connection and (\ref{Lie-2}).

Let $W$ be a vertical vector field on ${\cal Z}$. Then
$$
h_t(D_{V}X^h,W)=-h_t(X^h,D_{V}W)=0
$$
since the fibres are totally geodesic submanifolds, so $D_{V}W$ is a
vertical vector field. Therefore $D_{V}X^h$ is a horizontal vector
field. Moreover, $[V,X^h]$ is a vertical vector field, hence
$D_{V}X^h={\cal H}D_{X^h}V$. Thus
$$
h_t(D_{V}X^h,Y^h)=h_t(D_{X^h}V,Y^h)=-h_t(V,D_{X^h}Y^h).
$$
Now (\ref{D-vh}) follows from (\ref{D-hh}) and (\ref{r-r}).

\section{The second fundamental form of an almost Hermitian structure as a map into the twistor space}

Now let $J$ be an almost complex structure on the manifold $N$
compatible with the metric $g$. Then $J$ can be considered as a
section of the bundle $\pi:{\cal Z}\to N$. Thus we have a map
$J:(N,g)\to ({\cal Z},h_t)$ between Riemannian manifolds. Let
$J^{\ast}T{\cal Z}\to N$ be the pull-back of the bundle $T{\cal
Z}\to {\cal Z}$ under the map $J:N\to{\cal Z}$. Then we can consider
the differential $J_{\ast}:TN\to T{\cal Z}$ as a section of the
bundle $Hom(TN,J^{\ast}T{\cal Z})\to N$. Denote by $\widetilde D$
the connection on $J^{\ast}T{\cal Z}$ induced by the Levi-Civita
connection $D$ on $T{\cal Z}$. The Levi Civita connection $\nabla$
on $TN$ and the connection $\widetilde D$ on $J^{\ast}T{\cal Z}$
induce a connection $\widetilde\nabla$ on the bundle
$Hom(TN,J^{\ast}T{\cal Z})$. Recall that the second fundamental form
of the map $J$ is, by definition,
$$
\widetilde\nabla J_{\ast}.
$$
The map $J: (N,g)\to ({\cal Z},h_t)$ is harmonic if and only if
$$
\mathit{Trace}_{g}\widetilde\nabla J_{\ast}=0.
$$
Recall also that the map $J: (N,g)\to ({\cal Z},h_t)$ is totally
geodesic exactly when $\widetilde\nabla J_{\ast}=0$.

 Any (local) section $a$
of the bundle $A(TN)$ determines a (local) vertical vector field
$\widetilde a$ defined by
$$
{\widetilde a}_I=\frac{1}{2}(a(p)+I\circ a(p)\circ I),\quad
p=\pi(I).
$$
Thus if $aE_j=\sum\limits_{l=1}^n a_{jl}E_l$,
$$
{\widetilde a}=\textstyle{\sum\limits_{j<l}}\widetilde{a}_{jl}\frac{\partial}{\partial
y_{jl}}
$$
where
$$
\widetilde{a}_{jl}=\frac{1}{2}[a_{jl}\circ\pi+\textstyle{\sum\limits_{r,s=1}^n}
y_{jr}(a_{rs}\circ\pi)y_{sl}]
$$

The next lemma is "folklore".

\begin{lemma}\label{Xh-a til}
If $I\in{\cal Z}$ and $X$ is a vector field on a neighbourhood of
the point $p=\pi(I)$, then
$$
[X^h,\widetilde a]_I=(\widetilde{\nabla_{X}a})_I.
$$
\end{lemma}
{\bf Proof}. Take an orthonormal frame $E_1,...,E_n$ of $TN$ near
the point $p$ such that $\nabla E_i|_p=0$, $i=1,...,n$. Let $(\tilde
x_{i},y_{jl})$, $1\leq j<l\leq n$, be the local coordinates of
$A(TN)$ defined by means of a local coordinate system $x$ of $N$ at
$p$ and the frame $E_1,..,E_n$. Then, by (\ref{eq 3.1}),
$$
[X^h,\frac{\partial}{\partial y_{jl}}]_I=0,\quad j,l=1,...,n,
\quad X^{h}=\textstyle{\sum\limits_{i=1}^n} X^{i}(p)(\frac{\partial}{\partial\tilde
x_i})_I.
$$
It follows that
$$
[X^h,\widetilde a]_I=\frac{1}{2}[X_p(a_{jl})+\textstyle{\sum\limits_{k,m=1}^n}
y_{jk}(I)X_p(a_{km})y_{ml}(I)]=(\widetilde{\nabla_{X}a})_I
$$
since $$(\nabla_{X_p}a)(E_i)=\textstyle{\sum\limits_{l=1}^n} X_p(a_{jl})(E_l)_p.$$

\smallskip

\noindent {\bf Remark 1}. For every $I\in{\cal Z}$, we can find
local sections $a_1,...,a_{m^2-m}$ of $A(TN)$ whose values at
$p=\pi(I)$ constitute a basis of the vertical space ${\cal V}_I$ and
such that $\nabla a_{\alpha}|_p=0$, $\alpha=1,...,m^2-m$. Let
$\widetilde a_{\alpha}$ be the vertical vector fields determined by
the sections $a_{\alpha}$. Lemma~\ref{Xh-a til} and the Koszul
formula for the Levi-Civita connection imply that $h_t(D_{\widetilde
a_{\alpha}}\widetilde a_{\beta}, X^h)_I=0$ for every $X\in T_pN$.
Therefore, for every vertical vector fields $U$ and $V$, the
covariant derivative ($D_{U}V)_I$ at $I$ is a vertical vector. It
follows that the fibres of the twistor bundle are totally geodesic
submanifolds.

\smallskip

Let $I\in{\cal Z}$ and let $U,V\in{\cal V}_I$. Take sections $a$
and $b$ of $A(TN)$ such that $a(p)=U$, $b(p)=V$ for $p=\pi(I)$.
Let $\widetilde a$ and $\widetilde b$ be the vertical vector
fields determine by the sections $a$ and $b$. Taking into account
the fact that the fibre of ${\cal Z}$ through the point $I$ is a
totally geodesic submanifold  and applying formula (\ref{coder})
we get
\begin{equation}\label{a-b-til}
(D_{\widetilde a}\widetilde b)_I=\frac{1}{4}[UVI+IVU+I(UVI+IVU)I]=0.
\end{equation}

\begin{lemma}\label{v-frame}
For every $p\in N$, there exists a $h_t$-orthonormal frame of
vertical vector fields $\{V_{\alpha}:~\alpha=1,...,m^2-m\}$ such
that

\noindent $(1)$ $\quad (D_{V_{\alpha}}V_{\beta})_{J(p)}=0$,~~
$\alpha,\beta=1,...,m^2-m$.

\noindent $(2)$ $\quad$ If $X$ is a vector field near the point $p$,
$[X^h,V_{\alpha}]_{J(p)}=0$.

\noindent $(3)$ $\quad$  $\nabla_{X_p}(V_{\alpha}\circ J)\perp {\cal
V}_{J(p)}$

\end{lemma}

{\bf Proof}. Let $E_1,...,E_n$ be an orthonormal frame of $TN$ in a
neighbourhood $N$ of $p$ such that $J(E_{2k-1})_p=(E_{2k})_p$,
$k=1,...,m$, and $\nabla E_l|_p=0$, $l=1,...,n$. Define sections
$S_{ij}$,$1\leq i,j\leq n$ by (\ref{Sij}) and, as in
Section~\ref{prel}, set
$$
\begin{array}{c}
A_{r,s}=\frac{1}{\sqrt 2}(S_{2r-1,2s-1}-S_{2r,2s}),\quad
B_{r,s}=\frac{1}{\sqrt 2}(S_{2r-1,2s}+S_{2r,2s-1}),\\[6pt]
r=1,...,m-1,\>s=r+1,...,m.
\end{array}
$$
Then $\{(A_{r,s})_p,(B_{r,s})_p\}$ is a $G$-orthonormal basis of
the vertical space ${\cal V}_{J(p)}$ such that
$(B_{r,s})_p=J(A_{r,s})_p$. Note also that $\nabla
A_{r,s}|_p=\nabla B_{r,s}|_p=0$. Let $\widetilde A_{r,s}$ and
$\widetilde B_{r,s}$ be the vertical vector fields on ${\cal Z}$
determined by the sections $A_{r,s}$ and $B_{r,s}$ of $A(TN)$.
These vector fields constitute a frame of the vertical bundle
${\cal V}$ in a neighbourhood of the point $J(p)$.

Consider $\widetilde A_{r,s}\circ J$ as a section of $A(TN)$. Then,
if $X\in T_pN$, we have
$$
\begin{array}{c}
\nabla_{X_p}(\widetilde A_{r,s}\circ
J)=\frac{1}{2}\{(\nabla_{X_p}J)\circ
(A_{r,s})_p\circ J_p+J_p\circ (A_{r,s})\circ (\nabla_{X_p}J)\}\\[6pt]
=\frac{1}{2}\{-\nabla_{X_p}J\circ J_p\circ (A_{r,s})_p+J_p\circ (A_{r,s})\circ (\nabla_{X_p}J)\}\\[6pt]
=\frac{1}{2}[(B_{r,s})_p,\nabla_{X_p}J]
\end{array}
$$
The endomorphisms $(B_{r,s})_p$ and $\nabla_{X_p}J$ of $T_pN$ belong
to ${\cal V}_{J(p)}$,  so they anti-commute with $J(p)$, hence their
commutator commutes with $J(p)$. Therefore, in view of
(\ref{decom}), the commutator $[(B_{r,s})_p,\nabla_{X_p}J]$ is
$G$-orthogonal to the vertical space at $J(p)$. Thus
$$
\nabla_{X_p}(\widetilde A_{r,s}\circ J)\perp{\cal V}_{J(p)}
$$
and similarly  $\nabla_{X_p}(\widetilde B_{r,s}\circ J)\perp{\cal
V}_{J(p)}$.

It is convenient to denote the elements of the frame $\{\widetilde
A_{r,s},\widetilde B_{r,s}\}$ by $\{\widetilde V_1,...,\widetilde
V_{m^2-m}\}$. In this way we have a frame of vertical vector fields
near the point $J(p)$ with the property $(3)$ of the lemma.
Properties $(1)$ and $(2)$ are also satisfied by this frame
according to (\ref{a-b-til}) and Lemma~\ref{Xh-a til}, respectively.
In particular,
$$
(\widetilde V_{\gamma})_{J(p)}(h_t(\widetilde V_{\alpha},\widetilde
V_{\beta}))=0,\quad \alpha,\beta,\gamma=1,...,m^2-m.
$$

Note also that, in view of (\ref{D-vh}),
$$
{\cal V}(D_{X^h}\widetilde V_{\alpha})_{J(p)}=[X^h,\widetilde
V_{\alpha}]_{J(p)}=0,
$$
hence
$$
X^h_{J(p)}(h_t(\widetilde V_{\alpha},\widetilde V_{\beta}))=0.
$$

 Now it is clear that the $h_t$-orthonormal frame
$\{V_1,...,V_{m^2-m}\}$ obtained from $\{\widetilde
V_1,...,\widetilde V_{m^2-m}\}$ by the Gram-Schmidt process has the
properties stated in the lemma.

\begin{prop}\label{covder-dif}
For every  $X,Y\in T_pN$, $p\in N$,
$$
\begin{array}{c}
\widetilde\nabla J_{\ast}(X,Y) =\displaystyle{\frac{1}{2}}{\cal
V}(\nabla^{2}_{XY}J + \nabla^{2}_{YX}J)\\[8pt]
-\displaystyle{\frac{2t}{n}}[(R((J\circ
\nabla_{X}J)^{\wedge})Y)^h_{J(p)}
+(R((J\circ\nabla_{Y}J)^{\wedge})X)^h_{J(p)}]
\end{array}
$$
where $\nabla^{2}_{XY}J=\nabla_X\nabla_Y J-\nabla_{\nabla_XY}J$ is
the second covariant derivative of $J$.
\end{prop}

{\bf Proof}. Extend $X$ and $Y$ to vector fields in a
neighbourhood of the point $p$. Let $V_1,...,V_{m^2-m}$ be a
$h_t$-orthonormal frame of vertical vector fields with the
properties $(1)$ - $(3)$ stated in Lemma~\ref{v-frame}.

 We have
$$
J_{\ast}\circ Y=Y^h\circ J+\nabla_{Y}J=Y^h\circ
J+\textstyle{\sum\limits_{\alpha=1}^{m^2-m}}h_t(\nabla_{Y}J,V_\alpha\circ
J)(V_{\alpha}\circ J),
$$
hence
$$
\begin{array}{c}
{\widetilde D}_{X}(J_{\ast}\circ Y)=(D_{J_{\ast}X}Y^h)\circ J+
\sum\limits_{\alpha=1}^{m^2-m}h_t(\nabla_{Y}J,V_{\alpha})(D_{J_{\ast}X}V_{\alpha})\circ J\\[8pt]
+t\sum\limits_{\alpha=1}^{m^2-m}G(\nabla_X\nabla_{Y} J,V_{\alpha}\circ
J)(V_{\alpha}\circ J)
\end{array}
$$
This, in view of Lemma~\ref{LC}, implies
$$
\begin{array}{c}
\widetilde{D}_{X_p}(J_{\ast}\circ
Y)=(\nabla_{X}Y)^h_{J(p)}+\displaystyle{\frac{1}{2}}R(X\wedge
Y)J(p)\\[8pt]-\displaystyle{\frac{2t}{n}}(R((J\circ\nabla_{X}J)^{\wedge})Y)^h_{J(p)}
+t\sum\limits_{\alpha=1}^{m^2-m}G(\nabla_{X_p}\nabla_{Y}J,V_{\alpha}\circ
J)_{p}V_{\alpha}(J(p))\\[8pt]
-\displaystyle{\frac{2t}{n}}(R((J\circ\nabla_{Y}J)^{\wedge})X)^h_{J(p)}\\[8pt]
=(\nabla_{X_p}Y)^h_{J(p)} +\displaystyle{\frac{1}{2}}{\cal
V}(\nabla_{X_p}\nabla_{Y}J +\nabla_{Y_p}\nabla_{X}J)
+\frac{1}{2}\nabla_{[X,Y]_p}J \\[8pt]
-\displaystyle{\frac{2t}{n}}[R((J\circ\nabla_{X}J)^{\wedge})Y)^h_{J(p)}
+(R((J\circ\nabla_{Y}J)^{\wedge})X)^h_{J(p)}].
\end{array}
$$
It follows that
$$
\begin{array}{c}
\widetilde\nabla J_{\ast}(X,Y)=\widetilde D_{X_p}(J_{\ast}\circ Y)-
(\nabla_{X}Y)^h_{\sigma}-\nabla_{\nabla_{X_p}Y} J\\[8pt]
=\displaystyle{\frac{1}{2}}{\cal
V}(\nabla_{X_p}\nabla_{Y}J-\nabla_{\nabla_{X_p}Y} J
+\nabla_{Y_p}\nabla_{X}J--\nabla_{\nabla_{Y_p}X} J)\\[8pt]
-\displaystyle{\frac{2t}{n}}[R((J\circ\nabla_{X}J)^{\wedge})Y)^h_{J(p)}
+(R((J\circ\nabla_{Y}J)^{\wedge})X)^h_{J(p)}].
\end{array}
$$

\begin{cotmb}
If $(N,g,J)$ is K\"ahler, the map $J:(N,g)\to ({\cal Z},h_t)$ is a
totally geodesic isometric imbedding.
\end{cotmb}

\noindent {\bf Remark 2}. By a result of C. Wood \cite{W1,W2}, $J$
is a harmonic almost complex structure, i.e.  a harmonic section of
the twistor space $({\cal Z},h_t)\to (N,g)$ if and only if
$[J,\nabla^{\ast}\nabla J]=0$ where $\nabla^{\ast}\nabla$ is the
rough Laplacian. This, in view of the decomposition (\ref{decom}),
is equivalent to the condition that the vertical part of
$\nabla^{\ast}\nabla J=-\mathit{Trace}\nabla^2 J$ vanishes. Thus, by
Proposition~\ref{covder-dif}, $ J$ is a harmonic section if and only
if
$$
{\cal V} \mathit{Trace}\widetilde\nabla J_{\ast}=0.
$$

\section{The Atiyah-Hitchin-Singer and Eells-Salamon almost complex
structures as harmonic sections}

Let $(M,g)$ be an oriented Riemannian manifold of dimension four.
The twistor space of such a manifold has two connected components,
which can be identified with the unit sphere subbundles ${\cal
Z}_{\pm}$ of the bundles $\Lambda^2_{\pm}TM\to M$, the
eigensubbundles of the bundle $\pi:\Lambda^2TM\to M$ corresponding
to the eigenvalues $\pm 1$ of the Hodge star operator. The sections
of ${\cal Z}_{\pm}$ are the almost complex structures on $M$
compatible with the metric and $\pm$-orientation of $M$. The spaces
${\cal Z}_{+}$ and ${\cal Z}_{-}$ are  called the "positive" and the
"negative" twistor space of $(M,g)$.

The Levi-Civita connection $\nabla$ of $M$ preserves the bundles
$\Lambda^2_{\pm}TM$, so it induces a metric connection on each of
them denoted again by $\nabla$. The  horizontal distribution of
$\Lambda^2_{\pm}TM$ with respect to $\nabla$ is tangent to the
twistor space ${\cal Z}_{\pm}$. Thus we have the decomposition
$T{\cal Z}_{\pm}={\cal H}\oplus {\cal V}$ of the tangent bundle of
${\cal Z}_{\pm}$ into horizontal and vertical components. The
vertical space ${\cal V}_{\tau}=\{V\in T_{\tau}{\cal Z}_{\pm}:~
\pi_{\ast}V=0\}$ at a point $\tau\in{\cal Z}$ is the tangent space
to the fibre of ${\cal Z}_{\pm}$ through $\tau$. Considering
$T_{\tau}{\cal Z}_{\pm}$ as a subspace of
$T_{\tau}(\Lambda^2_{\pm}TM)$ (as we shall always do), ${\cal
V}_{\tau}$ is the orthogonal complement of $\tau$ in
$\Lambda^2_{\pm}T_{\pi(\tau)}M$.

Given $a\in \Lambda^2TM$, define, as in Sec.~\ref{CLCS}, an
endomorphism $K_a$ of $T_{\pi(a)}M$ by
$$
g(K_aX,Y)=2g(a,X\wedge Y),\quad X,Y\in T_{\pi(a)}M.
$$

For $\sigma\in{\cal Z}_{\pm}$, $K_{\sigma}$ is a complex structure
on the vector space $T_{\pi(\sigma)}M$ compatible with the metric
and $\pm$ the orientation.

Denote by $\times$ the vector-cross product in the $3$-dimensional
oriented Euclidean space $(\Lambda^2_{\pm}T_pM,g_p)$, $p\in M$.

It is easy to show that if $a,b\in\Lambda^2_{\pm}TM$
\begin{equation}\label{KaKb}
K_a\circ K_b=-g(a,b)Id \pm K_{a\times b}.
\end{equation}

This identity implies that for every vertical vector $V\in {\cal
V}_{\sigma}$ and every $X,Y\in T_{\pi(\sigma)}M$
\begin{equation}\label{aux}
g(V,X\wedge K_{\sigma}Y)=g(V,K_{\sigma}X\wedge Y)=g(\sigma\times
V,X\wedge Y).
\end{equation}

Note also that, in view of (\ref{calJ}), the Atiyah-Hitchin-Singer
and Eells-Salamon almost complex structures ${\cal J}_1$ and ${\cal
J}_2$ at a point $\sigma\in {\cal Z}_{\pm}$ can be written as
$$
\begin{array}{c}
{\cal J}_kV=\pm(-1)^{k+1}\sigma\times V~~\mbox{for}~~V\in{\cal
V}_{\sigma},\\[6pt]
{\cal J}_kX^h_{\sigma}=K_{\sigma}X~~\mbox{for}~~X\in
T_{\pi(\sigma)}M,\\[6pt]
 k=1,2.
\end{array}
$$

Denote by ${\cal B}:\Lambda^2TM\to \Lambda^2TM$ the endomorphism
corresponding to the traceless Ricci tensor. If  $s$ denotes the
scalar curvature of $(M,g)$ and $\rho:TM\to TM$ the Ricci operator,
$g(\rho(X),Y)=Ricci(X,Y)$, we have
$$
{\cal B}(X\wedge Y)=\rho(X)\wedge
Y+X\wedge\rho(Y)-\frac{s}{2}X\wedge Y.
$$
Let ${\cal W}: \Lambda^2TM\to \Lambda^2TM$ be the endomorphism
corresponding the Weyl conformal tensor. Denote the restriction of
${\cal W}$ to $\Lambda^2_{\pm}TM$ by ${\cal W}_{\pm}$, so ${\cal
W}_{\pm}$ sends $\Lambda^2_{\pm}TM$ to $\Lambda^2_{\pm}TM$ and
vanishes on $\Lambda^2_{\mp}TM$.

 It is well known that the curvature operator decomposes as (see
e.g. \cite[Chapter 1 H]{Besse})
$$
{\cal R}=\frac{s}{6}Id+{\cal B}+{\cal W}_{+}+{\cal W}_{-}.
$$
Note that this differ by the factor $1/2$ from \cite{Besse} because
of the factor $1/2$ in our definition of the induced metric on
$\Lambda^2TM$.

The Riemannian manifold $(M,g)$ is Einstein exactly when ${\cal
B}=0$. It is called self-dual (anti-self-dual) if ${\cal W}_{-}=0$
(resp. ${\cal W}_{+}=0$). By a well-known result of
Atiyah-Hitchin-Singer  \cite{AHS},  the almost complex structure
${\cal J}_1$ on ${\cal Z}_{-}$ (resp. ${\cal Z}_{+}$) is integrable
(i.e. comes from a complex structure) if and only if $(M,g)$ is
self-dual (resp. anti-self-dual). On the other hand the almost
complex structure ${\cal J}_2$ is never integrable by a result of
Eells-Salamon \cite{ES} but nevertheless it is very useful in
harmonic map theory.

\smallskip

{\noindent}{\bf Convention}. In what follows the negative twistor
space ${\cal Z}_{-}$ will be called simply "the twistor space" and
will be denoted by ${\cal Z}$.

Changing the orientation of $M$ interchanges the roles of
$\Lambda^2_{+}TM$ and $\Lambda^2_{-}TM$, respectively of ${\cal
Z}_{+}$ and ${\cal Z}_{-}$. But note that the Fubini-Study metric on
${\Bbb C}{\Bbb P}^2$ is self-dual and not anti-self-dual, so the
structure ${\cal J}_1$ on the negative twistor space ${\cal Z}_{-}$
is integrable while on ${\cal Z}_{+}$ it is not. This is one of the
reasons to prefer ${\cal Z}_{-}$ rather than ${\cal Z}_{+}$.

\smallskip

Remark 2, Proposition~\ref{covder-dif} and Theorem~1
imply
\begin{cotmb}\label{hs}

\noindent $(i)$ ${\cal V} \mathit{Trace}\widetilde\nabla {\cal
J}_{1\,\ast}=0$ if and only if $(M,g)$ is self-dual.

\noindent $(ii)$ ${\cal V}\mathit{Trace}\widetilde\nabla {\cal
J}_{2\,\ast}=0$ if and only if $(M,g)$ is self-dual and with
constant scalar curvature.
\end{cotmb}

\section{The Atiyah-Hitchin-Singer and Eells-Salamon almost complex
structures as harmonic maps}

In this section we prove Theorem 2, which is the main
result of the paper.

Note first that the almost complex structure ${\cal J}_k$, $k=1$ or $2$, is a
harmonic map if and only if ${\cal V}\mathit{Trace}\widetilde\nabla
{\cal J}_{k\,\ast}=0$ and ${\cal H}\mathit{Trace}\widetilde\nabla {\cal
J}_{k\,\ast}=0$. By Corollary~\ref{hs} if the vertical part of
$\mathit{Trace}\widetilde\nabla {\cal J}_{k\,\ast}$ vanishes, then the
manifold $(M,g)$ is self dual. According to
Proposition~\ref{covder-dif} ${\cal H}\mathit{Trace}\widetilde\nabla
{\cal J}_{k\,\ast}=0$, $k=1,2$, if and only if for every
$\sigma\in{\cal Z}$ and every $F\in T_{\sigma}{\cal Z}$
$$
\mathit{Trace}_{h_t}\,\{T_{\sigma}{\cal Z}\ni A\to h_t(R_{\cal Z}(({\cal
J}_{k}\circ D_{A}{\cal J}_k)^{\wedge})A),F)\}=0.
$$

Set for brevity
$$
Tr_k(F)=\mathit{Trace}_{h_t}\,\{T_{\sigma}{\cal Z}\ni A\to h_t(R_{\cal
Z}(({\cal J}_{k}\circ D_{A}{\cal J}_k)^{\wedge})A),F)\}.
$$

The next two technical lemmas, giving explicit formulas for
$Tr_k(F)$ in the self-dual case, will be proved in the next
section.

\begin{lemma}\label{tr-ver}
Suppose that $(M,g)$ is self-dual. Then, if $\sigma\in{\cal Z}$
and $U\in{\cal V}_{\sigma}$,
$$
Tr_k(U)=\frac{t}{4}g({\cal B}(U),{\cal B}(\sigma)), \> k=1,2.
$$
\end{lemma}

\begin{lemma}\label{tr-horr}
Suppose that $(M,g)$ is self-dual. Then, if $X\in T_{p}M$,
$p=\pi(\sigma)$,
$$
\begin{array}{c}
Tr_k(X^h_{\sigma})=[1+(-1)^k]\displaystyle{\frac{s(p)}{144}}X(s)
+\displaystyle{\frac{1}{12}}(\frac{ts(p)}{6}-2)X(s)\\[6pt]
+\mathit{Trace}_{h_t}\,\{{\cal V}_{\sigma}\ni V\to
[\displaystyle{\frac{t}{8}}g((\nabla_{X}{\cal B})( V),{\cal B}(V))\\[6pt]
\hspace{6.5cm}+(-1)^{k+1}\displaystyle{\frac{ts(p)}{24}}g(\delta{\cal
B}(K_VX),V)]\}.
\end{array}
$$
\end{lemma}

\smallskip

\noindent {\bf Proof of Theorem~2}.  Suppose that ${\cal
J}_1$ or ${\cal J}_2$ is a harmonic map. By Corollary~\ref{hs},
$(M,g)$ is self-dual or self-dual with constant  scalar curvature.
Moreover, $Tr_k(U)=0$ for every vertical vector $U$ and
$Tr_k(X^h)=0$ for every horizontal vector $X^h$, $k=1$ or $k=2$.
Note that in both cases the first term in the expression for
$Tr_k(X^h)$  given in Lemma~\ref{tr-horr} vanishes:
$$
[1+(-1)^k]\displaystyle{\frac{s(p)}{144}}X(s)=0,\quad k=1,2.
$$
By Lemma~\ref{tr-ver}, for every $p\in M$ and every orthonormal
basis $v_1,v_2,v_3$ of $\Lambda^2_{-}T_pM$, $g({\cal B}(v_i),{\cal
B}(v_j))=0$,  $i,j=1,2,3$, $i\neq j$. This implies $g({\cal
B}(v_i),{\cal B}(v_i))=g({\cal B}(v_j),{\cal B}(v_j))$, $i\neq j$.
It follows that the function ${\cal Z}_p\ni\sigma\to ||{\cal
B}(\sigma)||^2$ is constant on the fibre ${\cal Z}_p$ of ${\cal
Z}$ at $p$. Thus we have a smooth function $f$ on $M$ such that
$f(p)=||{\cal B}(\sigma)||^2$ for every $\sigma\in{\cal Z}_p$. It
follows that
\begin{equation}\label{Xf}
X(f)=2g((\nabla_{X}{\cal B})(\sigma),{\cal
B}(\sigma))
\end{equation}
for every tangent vector $X\in T_pM$.

Let $E_1,...,E_4$ be an oriented orthonormal basis of $T_pM$
consisting of eigenvectors of $\rho$. Denote by
$\lambda_1,...,\lambda_4$ the corresponding eigenvalues. We have
$\lambda_1+\lambda_2+\lambda_3+\lambda_4=s$ and
\begin{equation}\label{cal B}
{\cal B}(X\wedge Y)=\rho(X)\wedge
Y+X\wedge\rho(Y)-\frac{s}{2}X\wedge Y.
\end{equation}
Define $s_i^{+}$ and $s_i=s_i^{-}$, $i=1,2,3$, as in
(\ref{s-basis}) by means of the basis $E_1,...,E_4$. Then
$$
\begin{array}{c}
{\cal
B}(s_1)=(\lambda_1+\lambda_2-\displaystyle{\frac{s}{2}})s_1^{+},\quad
{\cal
B}(s_2)=(\lambda_1+\lambda_3-\displaystyle{\frac{s}{2}})s_2^{+},\\[6pt]
{\cal
B}(s_3)=(\lambda_1+\lambda_4-\displaystyle{\frac{s}{2}})s_3^{+}.
\end{array}
$$
Therefore $||{\cal B}(\cdot)||^2=const$ on the fibre ${\cal Z}_p$
if and only if
$$
|\lambda_1+\lambda_2-\displaystyle{\frac{s}{2}}|=|\lambda_1+\lambda_3-\displaystyle{\frac{s}{2}}|
=|\lambda_1+\lambda_4-\displaystyle{\frac{s}{2}}|,
$$
i.e. if and only if, at every point $p\in M$, three eigenvalues of
${\rho}$ coincide.

Moreover,
$$
3f(p)=||{\cal B}(s_1)||^2+||{\cal B}(s_2)||^2+||{\cal
B}(s_3)||^2=||\rho||^2-\frac{s^2(p)}{4}.
$$
and, by (\ref{Xf}), it follows  that
$$
\mathit{Trace}_{h_t}\,\{{\cal V}_{\sigma}\ni V\to g((\nabla_{X}{\cal
B})(V),{\cal B}(V))\}=\frac{1}{3}X(||\rho||^2)-\frac{s(p)X(s)}{6}.
$$
Fix a tangent vector $X\in T_pM$ and denote by  $P$  the symmetric
bilinear form on $\Lambda^2_{-}T_pM$ corresponding to the
quadratic form
\begin{equation}\label{P}
P(a,a)=\frac{ts(p)}{24}g(\delta{\cal B}(K_aX),a).
\end{equation}
Set
$$
\psi=-(\frac{ts(p)}{144}+\frac{1}{6})X(s)+\frac{t}{24}X(||\rho||^2).
$$
Then for every $\sigma\in{\cal Z}_p$ and every $V\in{\cal
V}_{\sigma}$ with $||V||_g=1$ we have
\begin{equation}\label{tr-last}
Tr_k(X^h_{\sigma})=(-1)^{k+1}[\frac{1}{t}P(V,V)+\frac{1}{t}P(\sigma\times
V,\sigma\times V)]+\psi.
\end{equation}
Let $\{s_1,s_2,s_3\}$ be an orthonormal basis of $\Lambda^2_{-}T_pM$.
Take
$$
\sigma=\frac{1}{\sqrt{y_1^2+y_2^2+y_3^2}}(y_1s_1+y_2s_2+y_3s_3)
$$
for $(y_1,y_2,y_3)\in{\Bbb R}^3$ with $y_1\neq 0$. Set
$$
V=\frac{1}{\sqrt{y_1^2+y_2^2}}(-y_2s_1+y_1s_2).
$$
Then
$$
\sigma\times
V=\frac{1}{\sqrt{(y_1^2+y_2^2)(y_1^2+y_2^2+y_3^2)}}(-y_1y_3s_1-y_2y_3s_2+(y_1^2+y_2^2)s_3).
$$
Now varying $(y_1,y_2,y_3)$ we see from (\ref{tr-last}) that the
identity $Tr_k(X^h_{\sigma})=0$ implies
$$
P(s_i,s_j)=0,\quad
(-1)^{k+1}\frac{1}{t}[P(s_i,s_i)+P(s_j,s_j)]+\psi=0,\quad
i,j=1,2,3,\>i\neq j.
$$
Since $P(s_i,s_j)=0$, $i\neq j$, for every orthonormal basis, we
have $P(s_i,s_i)=P(s_j,s_j)$.

Suppose that $s(p)\neq 0$. Then, by the latter identity,
$$
g(\delta{\cal B}(K_{s_i}X),s_i)=g(\delta{\cal
B}(K_{s_j}X),s_j),\quad i,j=1,2,3.
$$
Take an oriented orthonormal basis $E_1,...,E_4$ of $T_pM$ and,
using it, define $s_i=s_i^{-}$, $i=1,2,3$. Then $g(\delta{\cal
B}(K_{s_1}X),s_1)=g(\delta{\cal B}(K_{s_2}X),s_2)$ for every $X\in
T_pM$. This, in view of (\ref{KaKb}), gives
$$
-g(\delta{\cal B}(X),s_1)=g(\delta{\cal B}(K_{s_3}X),s_2),\quad
X\in T_pM.
$$
Applying the latter identity for the basis $E_3,E_4,E_1,E_2$ we
get
$$
g(\delta{\cal B}(X),s_1)=g(\delta{\cal B}(K_{s_3}X),s_2).
$$
Hence $g(\delta{\cal B}(X),s_1)=0$. Similarly $g(\delta{\cal
B}(X),s_2)=(\delta{\cal B}(X),s_3)=0$. Therefore for every $X\in
T_pM$ and $a\in\Lambda^2_{-}T_pM$
$$
g(\delta{\cal B}(X),a)=0.
$$
Then, by (\ref{P}), $P(a,a)=0$ for every $a\in\Lambda^2_{-}T_pM$.
Thus, we see from (\ref{tr-last}) that the condition
$Tr_k(X^h_{\sigma})=0$ for every $\sigma\in{\cal Z}$, $X\in
T_{\pi(\sigma)}M$ is equivalent to the identities
$$
g(\delta{\cal B}(X),\sigma)=0,\quad \psi=0.
$$

Identity (\ref{cal B}) implies that for every $X\in T_pM$ and
every orthonormal basis $E_1,..,E_4$ of $T_pM$
$$
\delta{\cal B}(X)=\delta\rho\wedge
X-\textstyle{\sum\limits_{m=1}^4}[E_m\wedge(\nabla_{E_m}\rho)(X)-\frac{1}{2}E_m(s)E_m\wedge
X].
$$
Therefore the identity $g(\delta{\cal B}(X),\sigma)=0$ is
equivalent to
$$
g(\delta\rho,K_{\sigma}X)+\textstyle{\sum\limits_{m=1}^4}g((\nabla_{E_m}\rho)(X),K_{\sigma}E_m)+\frac{1}{2}(K_{\sigma}X)(s)=0.
$$
This is equivalent to
\begin{equation}\label{newtr}
\textstyle{\sum\limits_{m=1}^4}g((\nabla_{E_m}\rho)(K_{\sigma}E_m),X)=0
\end{equation}
since $g(\delta\rho,Z)=-\frac{1}{2}Z(s)$ by the second Bianchi
identity and the Ricci operator $\rho$ is $g$-symmetric. Let
$r(X,Y)$ be the Ricci tensor and set
$$
dr(X,Y,Z)=(\nabla_{Y}r)(Z,X)-(\nabla_{Z}r)(Y,X).
$$
Thus
$$
dr(X,Y,Z)=g((\nabla_{Y}\rho)(Z),X)-g((\nabla_{Z}\rho)(Y),X).
$$
The left-hand side of (\ref{newtr}) clearly does not depend on the
choice of the basis $(E_1,...,E_4)$. So, take an oriented
orthonormal basis $(E_1,...,E_4)$ such that $E_2=K_{\sigma}E_1$
and $E_4=-K_{\sigma}E_3$. Then
$$
dr(X,E_1,E_2)-dr(X,E_3,E_4)=\textstyle{\sum\limits_{m=1}^4}g((\nabla_{E_m}\rho)(K_{\sigma}E_m),X).
$$
Denote by $W_{-}$ the $4$-tensor corresponding to the operator
${\cal W}_{-}$,
$$
W_{-}(X,Y,Z,T)=g({\cal W}_{-}(X\wedge Y),Z\wedge T).
$$
Then the second Bianchi identity implies
$$
dr(X,E_1,E_2)-dr(X,E_3,E_4)=-2[\delta W_{-}(X,E_1,E_2)-\delta
W_{-}(X,E_3,E_4)].
$$
Since $(M,g)$ is self-dual, we see from the latter identity that
identity (\ref{newtr}) is always satisfied. The above identity
shows also that
\begin{equation}\label{dr=0}
dr(X,\sigma)=0,\quad \sigma\in{\cal Z},\> X\in T_{\pi(\sigma)}M.
\end{equation}

Let $\lambda_1(p)\leq \lambda_2(p)\leq\lambda_3(p)\leq \lambda_4(p)$
be the eigenvalues of the symmetric operator $\rho_p:T_pM\to T_pM$
in the ascending order. It is well-known that the functions
$\lambda_1,...,,\lambda_4$ are continuous (see, e.g. \cite[Chapter
Two, \S 5.7 ]{K} or \cite[Chapter I, \S 3]{R}). We have seen that,
at every point of $M$, at least three eigenvalues of the operator
$\rho$ coincide. The set $U$ of points at which exactly three
eigenvalues coincide is open by the continuity of
$\lambda_1,...,\lambda_4$. For every $p\in U$ denote the simple
eigenvalue of $\rho$ by $\lambda(p)$ and the triple eigenvalue by
$\mu(p)$, so the spectrum of $\rho$ is $(\lambda,\mu,\mu,\mu)$ with
$\lambda(p)\neq \mu(p)$ for every $p\in U$. As is well-known, the
implicit function theorem implies that the function $\lambda$ is
smooth. It is also well-known that, in a neighbourhood of every
point $p$ of $U$, there is a (smooth) unit vector field $E_1$ which
is an eigenvector of $\rho$ corresponding to $\lambda$. (for a proof
see \cite[Chapter 9, Theorem 7]{L}). Fix $p\in U$ and choose local
vector fields $E_2,E_3,E_4$ such that $(E_1,E_2,E_3,E_4)$ is an
oriented orthonormal frame. Let $\alpha$ be the dual $1$-form to
$E_1$, $\alpha(X)=g(E_1,X)$. Then
$$
r(X,Y)=(\lambda-\mu)\alpha(X)\alpha(Y)-\mu g(X,Y)
$$
in a neighbourhood of $p$. Note that the function
$\mu=\frac{1}{3}(s-\lambda)$ is also smooth. Hence the identity
$\delta r=-\frac{1}{2}ds$ reads as
\begin{equation}\label{delta r}
\begin{array}{c}
-E_1(\lambda-\mu)\alpha(X)-X(\mu)+(\lambda-\mu)[\delta\alpha.\alpha(X)-(\nabla_{E_1}\alpha)(X)]\\[6pt]
=-\frac{1}{2}[X(\lambda)+3X(\mu)],\quad X\in TU.
\end{array}
\end{equation}

Let $s_i=s_i^{-}$, $i=1,2,3$, be defined by means of
$E_1,...,E_4$. Taking into account that
$(\nabla_{X}\alpha)(E_1)=0$, we easily see that the identities
$dr(E_k,s_1)=0$, $k=1,2,3,4$, give
\begin{equation}\label{dr-s}
\begin{array}{c}
(\lambda-\mu)[(\nabla_{E_1}\alpha)(E_2)-(\nabla_{E_3}\alpha)(E_4)+(\nabla_{E_4}\alpha)(E_3)]-E_2(\lambda)=0\\[6pt]
(\lambda-\mu)(\nabla_{E_2}\alpha)(E_2)-E_1(\mu)=0,\quad
(\lambda-\mu)(\nabla_{E_2}\alpha)(E_3)-E_4(\mu)=0\\[6pt]
(\lambda-\mu)(\nabla_{E_2}\alpha)(E_4)+E_3(\mu)=0.
\end{array}
\end{equation}
The identities obtained from the latter ones by cycle permutations
of $E_2,E_3,E_4$ also hold as a consequence of the identities
$dr(E_k,s_2)=0$ and $dr(E_k,s_3)=0$. Thus
\begin{equation}\label{jj}
(\lambda-\mu)(\nabla_{E_j}\alpha)(E_j)=E_1(\mu),\quad j=2,3,4.
\end{equation}
Hence
\begin{equation}\label{delta-alpha}
(\lambda-\mu)\delta\alpha=-3E_1(\mu)
\end{equation}
Moreover, we have
$$
(\nabla_{E_3}\alpha)(E_4)=E_2(\mu),\quad
(\nabla_{E_4}\alpha)(E_3)=-E_2(\mu)
$$
 and the first identity of (\ref{dr-s}) gives
$$
(\lambda-\mu)(\nabla_{E_1}\alpha)(E_2)=E_2(\lambda)+2E_2(\mu).
$$
On the other hand identity (\ref{delta r}) implies
$$
(\lambda-\mu)(\nabla_{E_1}\alpha)(E_2)=\frac{1}{2}E_2(\lambda)+\frac{1}{2}E_2(\mu).
$$
It follows that
$$
0=\frac{1}{2}E_2(\lambda)+\frac{3}{2}E_2(\mu)=\frac{1}{2}E_2(s),
$$
so $E_2(s)=0$. Similarly, $E_3(s)=0$ and $E_4(s)=0$. Identity
(\ref{delta r}) for $X=E_1$ together with (\ref{delta-alpha})
implies $0=E_1(\lambda)+3E_1(\mu)=E_1(s)$. It follows that the
scalar curvature $s$ is locally constant on $U$. Then identity
$\psi=0$ implies that $||\rho||^2$ is locally constant. Thus in a
neighbourhood of every point $p\in U$, we have $\lambda+3\mu=a$ and
$\lambda^2+3\mu^2=b^2$ where $a$ and $b$ are some constants. It
follows that
$$
\mu=\frac{1}{12}(3a\pm\sqrt{12b^2-3a^2}).
$$
Note that $12b^2-3a^2\neq 0$ since otherwise we would have
$\mu=\frac{1}{4}a$, hence $\lambda=a-3\mu=\frac{1}{4}a=\mu$, a
contradiction. Since $\mu$ is continuous, we see that $\mu$ is
constant, hence $\lambda$ is also constant. Then, by (\ref{jj}),
$(\nabla_{E_j}\alpha)(E_j)=0$ for $j=2,3,4$ and the first equation
of (\ref{dr-s}) gives $(\nabla_{E_1}\alpha)(E_2)=0$. Similarly
$(\nabla_{E_1}\alpha)(E_3)=(\nabla_{E_1}\alpha)(E_4)=0$. Thus
$(\nabla_{X}\alpha)(E_j)=0$ for every $X$ and $j=2,3,4$. This and
the obvious identity $(\nabla_{X}\alpha)(E_1)=0$ imply that the
$1$-form $\alpha$ is parallel. It follows that the restriction of
the Ricci tensor to $U$ is parallel.

In the interior of the closed set $M\setminus U$ the eigenvalues
of the Ricci tensor coincide, hence the metric $g$ is Einstein on
this open set. Therefore the scalar curvature $s$ is locally
constant on $Int\,(M\setminus U)$ and the Ricci tensor is parallel
on it. Thus the Ricci tensor is parallel on the open set $U\cup
Int\,(M\setminus U)=M\setminus bU$, where $bU$ stands for the
boundary of $U$. Since $M\setminus bU$ is dense in $M$ it follows
the Ricci tensor is parallel on $M$. This implies that the
eigenvalues $\lambda_1\leq...\leq \lambda_4$ of the Ricci tensor
are constant. Thus either $M$ is Einstein or exactly three of the
eigenvalues coincide. Since $(M,g)$ is self-dual, in the second
case the simple eigenvalue $\lambda$ vanishes by \cite[Lemma
1]{DGM}. Therefore $M$ is locally the product of an interval in
${\Bbb R}$ and a $3$-dimensional manifold of constant curvature.

Conversely, suppose that $(M,g)$ is self dual and either Einstein
or locally is the product of an interval and a manifold of
constant curvature. Then at least three of the eigenvalues of the
Ricci tensor coincide which, as we have seen, implies that
$||{\cal B}(\cdot)||^2=const$ on every  fibre of ${\cal Z}$. It
follows that $g({\cal B}(\sigma),{\cal B}({\tau}))=0$ for every
$\sigma,\tau\in{\cal Z}$ with $g(\sigma,\tau)=0$. Therefore
$Tr_k(U)=0$ for every vertical vector $U$, $k=1,2$, by
Lemma~\ref{tr-ver}. Moreover, $Tr_k(X^h)=0$ by Lemma~\ref{tr-horr}
since the scalar curvature is constant and $\nabla{\cal B}=0$.

\smallskip

\noindent {\bf Remark 3}. According to Theorems~1 and 2, the conditions under which  ${\mathcal J_1}$ or ${\mathcal J_2}$  is  a harmonic section or a harmonic map do not depend on the parameter $t$ of the metric $h_t$. Taking certain special values of $t$, we can obtain a metric $h_t$ with nice properties (cf., for example, \cite{DM91, DGM, M}).

\section{Proofs of Lemmas~\ref{tr-ver} and \ref{tr-horr}}

Denote by $R_{\cal Z}$ the curvature tensor of the Riemannian
manifold $({\cal Z},h_t)$.

\smallskip

Let $\Omega_{k,t}(A,B)=h_t({\cal J}_kA,B)$ be the fundamental
$2$-form of the almost Hermitian manifold $({\cal Z},h_t,{\cal
J}_k)$, $k=1,2$. Then, for $A,B,C\in T_{\sigma}{\cal Z}$,
$$
h_t({\cal J}_{k}\circ D_{A}{\cal J}_k)^{\wedge},B\wedge
C)=-\frac{1}{2}h_t((D_{A}{\cal J}_k)(B),{\cal
J}_kC)=-\frac{1}{2}(D_{A}\Omega_{k,t})(B,{\cal J}_kC).
$$
\begin{lemma}\label{D-Omega}{\rm (\cite{M})}
Let $\sigma\in{\cal Z}$ and $X,Y\in T_{\pi(\sigma)}M$, $V\in{\cal
V}_{\sigma}$. Then
$$
(D_{X^h_{\sigma}}\Omega_{k,t})(Y^h_{\sigma},V)=\frac{t}{2}[(-1)^kg({\cal
R}(V),X\wedge Y)-g({\cal R}(\sigma\times V),X\wedge K_{\sigma}Y)].
$$
$$
(D_{V}\Omega_{k,t})(X^h_{\sigma},Y^h_{\sigma})=\frac{t}{2}g({\cal
R}(\sigma\times V),X\wedge K_{\sigma}Y +K_{\sigma}X\wedge
Y)+2g(V,X\wedge Y).
$$
Moreover, $(D_{A}\Omega_{k,t})(B,C)=0$ when $A,B,C$ are three
horizontal vectors at $\sigma$ or at least two of them are vertical.
\end{lemma}

\begin{cotmb}\label{JDJ}
Let $\sigma\in{\cal Z}$, $X\in T_{\pi(\sigma)}M$, $U\in {\cal
V}_{\sigma}$. If $E_1,...,E_4$ is an orthonormal basis of
$T_{\pi(\sigma)}M$ and $V_1,V_2$ is a $h_t$-orthonormal basis of
${\cal V}_{\sigma}$,
$$
\begin{array}{l}
({\cal J}_k\circ D_{X^h_{\sigma}}{\cal
J}_k)^{\wedge}=-\displaystyle{\frac{1}{2}}\textstyle{\sum\limits_{i=1}^4\sum\limits_{l=1}^2}[g({\cal
R}(\sigma\times V_l),X\wedge E_i)\\[6pt]
\hspace{5.5cm} +(-1)^k g({\cal R}(V_l),X\wedge
K_{\sigma}E_i)](E_i^h)_{\sigma}\wedge V_l
\end{array}
$$
$$
\begin{array}{l}
({\cal J}_k\circ D_{U}{\cal J}_k)^{\wedge}=\textstyle{\sum\limits_{1\leq i<j\leq
4}}[\displaystyle{\frac{t}{2}}g({\cal R}(\sigma\times U),E_i\wedge
E_j-K_{\sigma}E_i\wedge K_{\sigma}E_j)\\[8pt]
\hspace{6cm}-2g(U,E_i\wedge K_{\sigma}E_j)](E_i^h)_{\sigma}\wedge
(E_j^h)_{\sigma}.
\end{array}
$$
\end{cotmb}

The sectional curvature of the Riemannian manifold $({\cal Z},h_t)$
can be computed in terms of the curvature of the base manifold $M$
by means of the following formula.

\begin{prop}\label{sec}{\rm (\cite{DM91})}
Let $E,F\in T_{\sigma}{\cal Z}$ and $X=\pi_{\ast}E$,
$Y=\pi_{\ast}F$, $V={\cal V}E$, $W={\cal V}F$. Then
$$
\begin{array}{c}
h_t(R_{\cal Z}(E,F)E,F)=g(R(X,Y)X,Y)\\[6pt]
-tg((\nabla_{X}{\cal R})(X\wedge Y),\sigma\times
W)+tg((\nabla_{Y}{\cal R})(X\wedge Y),\sigma\times V)\\[6pt]
-3tg({\cal R}(\sigma),X\wedge Y)g(\sigma\times V,W)\\[6pt]
-t^2g(R(\sigma\times V)X,R(\sigma\times
W)Y)+\displaystyle{\frac{t^2}{4}}||R(\sigma\times W)X+R(\sigma\times
V)Y||^2\\[6pt]
-\displaystyle{\frac{3t}{4}}||R(X,Y)\sigma||^2+t(||V||^2||W||^2-g(V,W)^2).

\end{array}
$$
\end{prop}

Using this formula, the well-known expression of the Levi-Civita
curvature tensor by means of sectional curvatures and differential
Bianchi identity one gets the following.

\begin{cotmb}\label{RZ}
Let $\sigma\in{\cal Z}$, $X,Y,Z,T\in T_{\pi(\sigma)}M$, and
$U,V,W\in{\cal V}_{\sigma}$. Then
$$
\begin{array}{c}
h_t(R_{\cal
Z}(X^h,Y^h)Z^h,T^h)_{\sigma}=g(R(X,Y)Z,T)\\[6pt]
-\displaystyle\frac{3t}{12}[2g(R(X,Y)\sigma,R(Z,T)\sigma)
-g(R(X,T)\sigma,R(Y,Z)\sigma)\\[6pt]
+g(R(X,Z)\sigma,R(Y,T)\sigma)].
\end{array}
$$
$$
\begin{array}{c}
h_t(R_{\cal
Z}(X^h,Y^h)Z^h,U)_{\sigma}=-\displaystyle{\frac{t}{2}g(\nabla_{Z}{\cal
R}(X\wedge Y),\sigma\times U)}.\\[8pt]
h_t(R_{\cal
Z}(X^h,U)Y^h,V)_{\sigma}=\displaystyle{\frac{t^2}{4}g(R(\sigma\times
V)X,R(\sigma\times U)Y)}\\[6pt]
+\displaystyle{\frac{t}{2}g({\cal R}(\sigma),X\wedge
Y)g(\sigma\times V,U)}.\\[8pt]
h_t(R_{\cal Z}(X^h,Y^h)U,V)_{\sigma}
=\displaystyle{\frac{t^2}{4}}[g(R(\sigma\times V)X,R(\sigma\times
U)Y)\\[6pt]
\hspace{6cm}-g(R(\sigma\times U)X,R(\sigma\times V)Y)]\\[6pt]
+tg({\cal R}(\sigma),X\wedge Y)g(\sigma\times V,U).
\end{array}
$$
$$
h_t(R_{\cal Z}(X^h,U)V,W)=0.
$$
\end{cotmb}

We have stated in Lemma~\ref{tr-ver} that if $(M,g)$ is self-dual,
$$
Tr_k(U)=\frac{t}{4}g({\cal B}(U),{\cal B}(\sigma))~\rm{for}~\rm{
every}~ U\in{\cal V}_{\sigma},\> \sigma\in{\cal Z}.
$$

\noindent {\bf Proof of Lemma~\ref{tr-ver}}. Let $E_1,...,E_4$ be
an orthonormal basis of $T_pM$, $p=\pi(\sigma)$,  such that
$E_2=K_{\sigma}E_1$, $E_4=-K_{\sigma}E_3$. Define
$s_1=s_1^{-},s_2=s_2^{-},s_3=s_3^{-}$ via (\ref{s-basis}) by means
of $E_1,...,E_4$, so that $\sigma=s_1$ and ${\cal
V}_{\sigma}=span\{s_2,s_3\}$. Thus $V_1=\frac{1}{\sqrt t}s_2$,
$V_2=\frac{1}{\sqrt t}s_3$ is a $h_t$-orthonormal basis of ${\cal
V}_{\sigma}$.

\smallskip

By Corollary~\ref{JDJ}, for every $U\in {\cal V}_{\sigma}$
\begin{equation}\label{TrU}
\begin{array}{c}
Tr_k(U)=-\displaystyle{\frac{1}{2}}\textstyle{\sum\limits_{i,j=1}^4\sum\limits_{l=1}^2}[g({\cal
R}(\sigma\times V_l),E_j\wedge E_i)+(-1)^kg({\cal R}(V_l),E_j\wedge
K_{\sigma}E_i]\\[6pt]
\times h_t(R_{\cal Z}(E_i^h,V_l)E_j^h,U)\\[6pt]
+\displaystyle{\frac{t}{2}}\textstyle{\sum\limits_{l=1}^2}\{g({\cal R}(\sigma\times
V_l),s_2)[h_t(R_{\cal Z}(E_1^h,E_3^h)V_l,U)-h_t(R_{\cal
Z}(E_4^h,E_2^h)V_l,U)]\\[6pt]
\hspace{1.5cm}+g({\cal R}(\sigma\times V_l),s_3)[h_t(R_{\cal
Z}(E_1^h,E_4^h)V_l,U)-h_t(R_{\cal Z}(E_2^h,E_3^h)V_l,U)]\}\\[6pt]
-2\textstyle{\sum\limits_{1\leq i<j\leq 4}\sum\limits_{l=1}^2}g(V_l,E_i\wedge
K_{\sigma}E_j)h_t(R_{\cal Z}(E_i^h,E_j^h)V_l,U)
\end{array}
\end{equation}

We show first that
\begin{equation}\label{Isum}
\begin{array}{c}
\textstyle{\sum\limits_{i,j=1}^4\sum\limits_{l=1}^2}[g({\cal R}(\sigma\times V_l),E_j\wedge
E_i)h_t(R_{\cal Z}(E_i^h,V_l)E_j^h,U)\\[8pt]
=-\displaystyle{\frac{t}{2}}\mathit{Trace}_{h_t}\{{\cal V}_{\sigma}\ni V\to
g({\cal R}(\sigma\times V),{\cal R}(\sigma))g(\sigma\times U,V)\}.
\end{array}
\end{equation}

In order to prove this identity, we note that if
 $F\in T_{\sigma}{\cal Z}$, $V\in{\cal V}_{\sigma}$ and $a\in
\Lambda^2 T_{\pi(\sigma)}M$, the algebraic Bianchi identity
implies
\begin{equation}\label{a}
\begin{array}{c}
\textstyle{\sum\limits_{i,j=1}^4} g(a,E_j\wedge E_i)h_t(R_{\cal
Z}(E_i^h,V)E_j^h,F)\\[8pt]
=-\displaystyle{\frac{1}{2}}\textstyle{\sum\limits_{i,j=1}^4}g(a,E_i\wedge
E_j)h_t(R_{\cal Z}(E_i^h,E_j^h)V,F).
\end{array}
\end{equation}
Using the latter identity and  Corollary~\ref{RZ} we obtain
$$
\begin{array}{c}
\sum\limits_{i,j=1}^4\sum\limits_{l=1}^2[g({\cal R}(\sigma\times V_l),E_j\wedge
E_i)h_t(R_{\cal Z}(E_i^h,V_l)E_j^h,U)\\[8pt]
=-\displaystyle{\frac{t}{2}}\textstyle{\sum\limits_{l=1}^2}g({\cal R}(\sigma\times
V_l),{\cal R}(\sigma))g(\sigma\times U,V_l)\\[8pt]
=-\displaystyle{\frac{t}{2}}\mathit{Trace}_{h_t}\{{\cal V}_{\sigma}\ni V\to
g({\cal R}(\sigma\times V),{\cal R}(\sigma))g(\sigma\times U,V)\}.
\end{array}
$$
Next, we claim that
\begin{equation}\label{IIsum}
\textstyle{\sum\limits_{i,j=1}^4\sum\limits_{l=1}^2} g({\cal R}(V_l),E_j\wedge
K_{\sigma}E_i) h_t(R_{\cal Z}(E_i^h,V_l)E_j^h,U)=0.
\end{equation}
For every $V\in{\cal V}_{\sigma}$, we have
\begin{equation}\label{2dsum}
\begin{array}{c}
\sum\limits_{i,j=1}^4g({\cal R}(V),E_j\wedge K_{\sigma}E_i)h_t(R_{\cal
Z}(E_i^h,V)E_j^h,U)\\[8pt]
=g({\cal R}(V),E_1\wedge E_2)[h_t(R_{\cal
Z}(E_1^h,V)E_1^h,U)+h_t(R_{\cal Z}(E_2^h,V)E_2^h,U)]\\[6pt]
-g({\cal R}(V),E_3\wedge E_4)[h_t(R_{\cal
Z}(E_3^h,V)E_3^h,U)+h_t(R_{\cal Z}(E_4^h,V)E_4^h,U)]\\[6pt]
+g({\cal R}(V),E_1\wedge E_3)[h_t(R_{\cal
Z}(E_4^h,V)E_1^h,U)+h_t(R_{\cal Z}(E_2^h,V)E_3^h,U)]\\[6pt]
+g({\cal R}(V),E_1\wedge E_4)[-h_t(R_{\cal
Z}(E_3^h,V)E_1^h,U)+h_t(R_{\cal Z}(E_2^h,V)E_4^h,U)]\\[6pt]
+g({\cal R}(V),E_2\wedge E_3)[h_t(R_{\cal
Z}(E_4^h,V)E_2^h,U)-h_t(R_{\cal Z}(E_1^h,V)E_3^h,U)]\\[6pt]
+g({\cal R}(V),E_4\wedge E_2)[h_t(R_{\cal
Z}(E_3^h,V)E_2^h,U)+h_t(R_{\cal Z}(E_1^h,V)E_4^h,U)]\\[6pt]
\end{array}
\end{equation}
Corollary~\ref{RZ} implies that
$$
\begin{array}{c}
h_t(R_{\cal Z}(E_4^h,V)E_1^h,U)+h_t(R_{\cal
Z}(E_2^h,V)E_3^h,U)\\[8pt]
=\displaystyle{\frac{t^2}{4}}[g(R(\sigma\times U)E_4,E_2)g({\cal
R}(\sigma\times V),s_1)\\[6pt]
\hspace{5cm} +g(R(\sigma\times V)E_1,E_3)g({\cal R}(\sigma\times
U),s_1)]\\[6pt]
-\displaystyle{\frac{t}{2}}g({\cal R}(\sigma),s_3)g(\sigma\times
U,V)
\end{array}
$$
Since $(M,g)$ is self-dual, for every
$\tau\in\Lambda^2_{-}T_{\pi(\sigma)}M$,
$$
{\cal R}(\tau)=\frac{s}{6}\tau +{\cal B}(\tau)
$$
where ${\cal B}(\tau)\in \Lambda^2_{+}T_{\pi(\sigma)}M$. Therefore
$$
g({\cal R}(\sigma\times V),s_1)=g({\cal R}(\sigma\times V),\sigma)=0
$$
and
$$ g({\cal R}(\sigma\times U),s_1)=0,\quad g({\cal
R}(\sigma),s_3)=0.
$$
Thus
\begin{equation}\label{sim}
h_t(R_{\cal Z}(E_4^h,V)E_1^h,U)+h_t(R_{\cal Z}(E_2^h,V)E_3^h,U)=0.
\end{equation}
Similarly
\begin{equation}\label{3ple}
\begin{array}{l}
-h_t(R_{\cal Z}(E_3^h,V)E_1^h,U)+h_t(R_{\cal
Z}(E_2^h,V)E_4^h,U)=0\\[8pt]
h_t(R_{\cal Z}(E_4^h,V)E_2^h,U)-h_t(R_{\cal
Z}(E_1^h,V)E_3^h,U)=0\\[8pt]
h_t(R_{\cal Z}(E_3^h,V)E_2^h,U)+h_t(R_{\cal Z}(E_1^h,V)E_4^h,U)=0.
\end{array}
\end{equation}
Moreover, a straightforward computation gives
$$
\begin{array}{c}
\textstyle{\sum\limits_{l=1}^2 }\{g({\cal R}(V_l),E_1\wedge E_2)[h_t(R_{\cal
Z}(E_1^h,V_l)E_1^h,U)+h_t(R_{\cal Z}(E_2^h,V_l)E_2^h,U)]\\[6pt]
-g({\cal R}(V_l),E_3\wedge E_4)[h_t(R_{\cal
Z}(E_3^h,V_l)E_3^h,U)+h_t(R_{\cal Z}(E_4^h,V_l)E_4^h,U)]\}\\[8pt]
=\displaystyle{\frac{t^2}{8}}\textstyle{\sum\limits_{l=1}^2}g({\cal
R}(V_l),s_1^{+})g({\cal R}(\sigma\times U),s_1)g({\cal
B}(\sigma\times V_l),s_1^{+})\\[6pt]
=0.
\end{array}
$$
In view of (\ref{2dsum}), the latter identity, (\ref{sim}) and
(\ref{3ple}) imply (\ref{IIsum}).

\smallskip

Using the algebraic Bianchi identity, we see from (\ref{3ple})
that
$$
\begin{array}{c}
h_t(R_{\cal Z}(E_1^h,E_3^h)V,U)-h_t(R_{\cal
Z}(E_4^h,E_2^h)V,U)=0\\[6pt]
h_t(R_{\cal Z}(E_1^h,E_4^h)V,U)-h_t(R_{\cal Z}(E_2^h,E_3^h)V,U)=0.
\end{array}
$$
Hence
\begin{equation}\label{IIIsum}
\begin{array}{c}
\sum\limits_{l=1}^2\{g({\cal R}(\sigma\times V_l),s_2)[h_t(R_{\cal
Z}(E_1^h,E_3^h)V_l,U)-h_t(R_{\cal
Z}(E_4^h,E_2^h)V_l,U)]\\[6pt]
\hspace{1cm}+g({\cal R}(\sigma\times V_l),s_3)[h_t(R_{\cal
Z}(E_1^h,E_4^h)V_l,U)-h_t(R_{\cal Z}(E_2^h,E_3^h)V_l,U)]\}\\[6pt]
=0.
\end{array}
\end{equation}
Using (\ref{aux}) and Corollary~\ref{RZ}, we get
$$
\begin{array}{c}
\sum\limits_{1\leq i<j\leq 4}g(V,E_i\wedge
K_{\sigma}E_j)h_t(R_{\cal Z}(E_i^h,E_j^h)V,U)\\[8pt]
=\displaystyle{\frac{t^2}{4}}\textstyle{\sum\limits_{1\leq i<j\leq 4}}g(\sigma\times
V,E_i\wedge E_j)[g(R(\sigma\times U)E_i,R(\sigma\times
V)E_j)\\[6pt]
\hspace{7cm}-g(R(\sigma\times V)E_i,R(\sigma\times U)E_j]\\[9pt]
=\displaystyle{\frac{t^2}{8}}\textstyle{\sum\limits_{i=1}^4}g(R(\sigma\times
U)E_i,R(\sigma\times V)K_{\sigma\times V}E_i)
\end{array}
$$
Therefore
$$
\begin{array}{c}
\textstyle{\sum\limits_{1\leq i<j\leq 4}\sum\limits_{l=1}^2}g(V_l,E_i\wedge
K_{\sigma}E_j)h_t(R_{\cal Z}(E_i^h,E_j^h)V_l,U)\\[8pt]
=\displaystyle{\frac{t}{8}}\textstyle{\sum\limits_{i,k=1}^4}[g(R(\sigma\times
U)E_i,E_k)g(R(s_3)K_{s_3}E_i,E_k)\\[6pt]
\hspace{5cm}+g(R(\sigma\times U)E_i,E_k)g(R(s_2)K_{s_2}E_i,E_k)]\\[8pt]
=\displaystyle{\frac{t}{8}}[-g({\cal R}(\sigma\times
U),\sigma)g({\cal R}(s_3),s_2)+g({\cal R}(\sigma\times
U),s_2)g({\cal R}(s_3),s_1)\\[6pt]
\hspace{1.2cm}+g({\cal R}(\sigma\times U),\sigma)g({\cal
R}(s_2),s_3)-g({\cal R}(\sigma\times U),s_3)g({\cal
R}(s_2),s_1)]\\[8pt]
\end{array}
$$
This, by virtue of the self-duality of $(M,g)$, gives
\begin{equation}\label{IVsum}
\textstyle{\sum\limits_{1\leq i<j\leq 4}\sum\limits_{l=1}^2}g(V_l,E_i\wedge
K_{\sigma}E_j)h_t(R_{\cal Z}(E_i^h,E_j^h)V_l,U)=0.
\end{equation}
Identities (\ref{TrU}),(\ref{Isum}),(\ref{IIsum}), (\ref{IIIsum})
and (\ref{IVsum}) imply
$$
Tr_k(U)=\frac{t}{4}\mathit{Trace}_{h_t}\{{\cal V}_{\sigma}\ni V\to g({\cal
R}(\sigma\times V),{\cal R}(\sigma))g(\sigma\times U,V)\},\>
k=1,2.
$$
Now the lemma follows from the latter identity since $g({\cal
R}(\tau),{\cal R}(\sigma))=g({\cal B}(\tau),{\cal B}(\sigma))$ for
every $\tau, \sigma$ with $\tau\perp\sigma$.

\medskip

Recall that, according to Lemma~\ref{tr-horr}, if $(M,g)$ is
self-dual
$$
\begin{array}{c}
Tr_k(X^h_{\sigma})=[1+(-1)^k]\displaystyle{\frac{s(p)}{144}}X(s)
+\displaystyle{\frac{1}{12}}(\frac{ts(p)}{6}-2)X(s)\\[6pt]
+\mathit{Trace}_{h_t}\,\{{\cal V}_{\sigma}\ni V\to
[\displaystyle{\frac{t}{8}}g((\nabla_{X}{\cal B})( V),{\cal
B}(V))\\[6pt]
\hspace{6.5cm}+(-1)^{k+1}\displaystyle{\frac{ts(p)}{24}}g(\delta{\cal
B}(K_VX),V)]\}.
\end{array}
$$
for $X\in T_{\pi(\sigma)}$, $\sigma\in{\cal Z}$.

\smallskip

\noindent {\bf Proof of Lemma~\ref{tr-horr}}. Let
$s_1=s_1^{-},s_2=s_2^{-},s_3=s_3^{-}$ be the basis of
$\Lambda^2_{-}T_pM$, $p=\pi(\sigma)$, defined by means of an
oriented orthonormal basis $E_1,...,E_4$  of $T_pM$ such that
$E_2=K_{\sigma}E_1$, $E_4=-K_{\sigma}E_3$. Set $V_1=\frac{1}{\sqrt
t}s_2$, $V_2=\frac{1}{\sqrt t}s_3$.

Then, by Corollary~\ref{JDJ},
$$
\begin{array}{c}
Tr_k(X^h_{\sigma})=-\displaystyle{\frac{1}{2}}\textstyle{\sum\limits_{i,j=1}^4\sum\limits_{l=1}^2}[g({\cal
R}(\sigma\times V_l),E_j\wedge E_i)+(-1)^kg({\cal R}(V_l),E_j\wedge
K_{\sigma}E_i] \\[6pt]
\times h_t(R_{\cal Z}(E_i^h,V_l)E_j^h,X^h)\\[6pt]
+\sum\limits_{i<j}\sum\limits_{l=1}^2[\displaystyle{\frac{t}{2}}g({\cal
R}(\sigma\times V_l),E_i\wedge E_j-K_{\sigma}E_i\wedge
K_{\sigma}E_j)-2g(V_l,E_i\wedge K_{\sigma}E_j)]\\[6pt]
\times h_t(R_{\cal Z}(E_i^h,E_j^h)V_l,X^h).
\end{array}
$$

Identity (\ref{a}) and Corollary~\ref{RZ} imply
$$
\begin{array}{c}
\sum\limits_{i,j=1}^4\sum\limits_{l=1}^2g({\cal R}(\sigma\times V_l),E_j\wedge
E_i)h_t(R_{\cal Z}(E_i^h,V_l)E_j^h,X^h)_{\sigma}\\[6pt]
=-\displaystyle{\frac{t}{4}}\textstyle{\sum\limits_{l=1}^2}g(\displaystyle{\frac{1}{6}}X(s)\sigma\times
V_l+(\nabla_{X}{\cal B})(\sigma\times
V_l),\displaystyle{\frac{1}{6}}s(p)\sigma\times V_l+{\cal
B})(\sigma\times V_l)\\[6pt]
=-\displaystyle{\frac{s(p)}{72}}X(s)-\displaystyle{\frac{t}{4}}\mathit{Trace}_{h_t}\{{\cal
V}_{\sigma}\ni V\to g((\nabla_{X}{\cal B})(V),{\cal B}(V))\},
\end{array}
$$
where the latter identity follows from the fact that
$g((\nabla_{X}{\cal B})(a),b)=0$ for every
$a,b\in\Lambda^2_{-}T_pM$ (since the operator ${\cal B}$ sends
$\Lambda^2_{-}TM$ into $\Lambda^2_{+}TM$, and the connection
$\nabla$ preserves the bundles $\Lambda^2_{\pm}TM$).

Taking into account identity (\ref{aux}) and the fact that
$$
E_i\wedge E_j-K_{\sigma}E_i\wedge
K_{\sigma}E_j\in\Lambda^2_{-}T_{\pi(\sigma)}M,
$$
we have
$$
\begin{array}{c}
\sum\limits_{i<j}\sum\limits_{l=1}^2[\displaystyle{\frac{t}{2}}g({\cal
R}(\sigma\times V_l),E_i\wedge E_j-K_{\sigma}E_i\wedge
K_{\sigma}E_j)-2g(V_l,E_i\wedge K_{\sigma}E_j)]\\[6pt]
\times h_t(R_{\cal Z}(E_i^h,E_j^h)V_l,X^h)_{\sigma}\\[8pt]
=(\displaystyle{\frac{ts(p)}{6}}-2)\textstyle{\sum\limits_{i<j}\sum\limits_{l=1}^2}g(\sigma\times
V_l,E_i\wedge E_j)h_t(R_{\cal Z}(E_i^h,E_j^h)V_l,X^h)\\[8pt]
=\displaystyle{\frac{t}{4}(\frac{ts(p)}{6}-2)}\textstyle{\sum\limits_{l=1}^2}g((\nabla_{X}{\cal
R})(\sigma\times V_l),\sigma\times V_l)\\
=\displaystyle{\frac{1}{4}(\frac{ts(p)}{6}-2)\frac{X(s)}{3}}.
\end{array}
$$
Thus
$$
\begin{array}{c}
Tr_k(X^h_{\sigma})=(-1)^{k+1}\displaystyle{\frac{1}{2}}\textstyle{\sum\limits_{i,j=1}^4\sum\limits_{l=1}^2}g({\cal
R}(V_l),E_j\wedge K_{\sigma}E_ih_t(R_{\cal
Z}(E_i^h,V_l)E_j^h,X^h)_{\sigma}\\[8pt]
+\displaystyle{\frac{s(p)}{144}}X(s)
+\displaystyle{\frac{1}{12}(\frac{ts(p)}{6}}-2)X(s)\\[6pt]
+\mathit{Trace}_{h_t}\,\{{\cal V}_{\sigma}\ni V\to
\displaystyle{\frac{t}{8}}g((\nabla_{X}{\cal B})(V),{\cal B}(V))\}\\[6pt]
\end{array}
$$
In order to compute the first summand in the right-hand side of
the latter identity, it is convenient to set $C_{ilj}=h_t(R_{\cal
Z}(E_i^h,V_l)E_j^h,X^h)_{\sigma}$. Then
$$\begin{array}{c}
\textstyle{\sum\limits_{i,j=1}^4\sum\limits_{l=1}}g({\cal R}(V_l),E_j\wedge
K_{\sigma}E_i)h_t(R_{\cal
Z}(E_i^h,V_l)E_j^h,X^h)_{\sigma}\\[6pt]
=\displaystyle{\frac{1}{2}}\textstyle{\sum\limits_{l=1}^2}[g({\cal
R}(V_l),s_1^{+}+s_1)C_{1l1}-g({\cal R}(V_l),s_3^{+}-s_3)C_{1l3} +g({\cal R}(V_l),s_2^{+}-s_2)C_{1l4}\\[10pt]

+g({\cal R}(V_l),s_1^{+}+s_1)C_{2l2}+g({\cal R}(V),s_2^{+}+s_2)C_{2l3}+g({\cal R}(V),s_3^{+}+s_3)C_{2l4}\\[10pt]

-g({\cal R}(V),s_3^{+}+s_3)C_{3l1}+g({\cal R}(V),s_2^{+}-s_2)C_{3l2} -g({\cal R}(V),s_1^{+}-s_1)C_{3l3}\\[10pt]

g({\cal R}(V),s_2^{+}+s_2)C_{4l1}+g({\cal R}(V),s_3^{+}-s_3)C_{4l2}
-g({\cal R}(V),s_1^{+}-s_1)C_{4l4}]\\[6pt]
=\displaystyle{\frac{s(p)}{12\sqrt
t}}[(-C_{114}+C_{213}-C_{312}+C_{411})
+(C_{123}+C_{224}-C_{321}-C_{422})]\\[10pt]

+\displaystyle{\frac{1}2}\textstyle{\sum\limits_{l=1}^2} [g({\cal
B}(V_l),s_1^{+})(C_{1l1}+C_{2l2}-C_{3l3}-C_{4l4})\\[10pt]

\hspace{1.1cm}+g({\cal B}(V_l),s_2^{+})(C_{1l4}+C_{2l3}+C_{3l2}+C_{4l1})\\[10pt]

\hspace{1.5cm}+g({\cal
B}(V_l),s_3^{+})(-C_{1l3}+C_{2l4}-C_{3l1}+C_{4l2})].
\end{array}
$$
By Corollary~\ref{RZ}
$$
\begin{array}{c}
-C_{124}+C_{223}-C_{322}+C_{421}\\[6pt]
=\displaystyle{\frac{\sqrt
t}{2}\textstyle{\sum\limits_{i=1}^4}[-\frac{1}{12}E_i(s)g(E_i,X)+g((\nabla_{E_i}{\cal
B})(K_{s_3}E_i\wedge X),s_3)]}.
\end{array}
$$
For every $i=1,...,4$, $K_{s_3}E_i\wedge X+E_i\wedge
K_{s_3}X\in\Lambda^2_{-}T_pM$.  Hence
$$
g((\nabla_{E_i}{\cal B})(K_{s_3}E_i\wedge X+E_i\wedge
K_{s_3}X),s_3)=0.
$$
It follows that
$$
-C_{124}+C_{223}-C_{322}+C_{421}=\displaystyle{\frac{\sqrt
t}{2}[-\frac{1}{12}X(s)+g(\delta{\cal B}(K_{s_3}X),s_3)]}.
$$
Similarly
$$
\begin{array}{c}
C_{133}+C_{234}-C_{331}-C_{432}\\[6pt]
=\displaystyle{\frac{\sqrt t}{2}[-\frac{1}{12}X(s)+g(\delta{\cal
B}(K_{s_2}X),s_2)]}.
\end{array}
$$
Hence
$$
\begin{array}{c}
(-C_{124}+C_{223}-C_{322}+C_{421})
+(C_{133}+C_{234}-C_{331}-C_{432})\\[8pt]
=\displaystyle{\frac{\sqrt
t}{2}}[-\frac{1}{6}X(s)+t\mathit{Trace}_{h_t}\,\{{\cal V}_{\sigma}\ni V\to
g(\delta{\cal B}(K_VX),V)\}.
\end{array}
$$
Set for short
$$
\begin{array}{c}
\Sigma(E_1,...,E_4)=\textstyle{\sum\limits_{l=1}^2} [g({\cal
B}(V_l),s_1^{+})(C_{1l1}+C_{2l2}-C_{3l3}-C_{4l4})\\[10pt]

\hspace{0.8cm}+g({\cal B}(V_l),s_2^{+})(C_{1l4}+C_{2l3}+C_{3l2}+C_{4l1})\\[10pt]

\hspace{1.2cm}+g({\cal
B}(V_l),s_3^{+})(-C_{1l3}+C_{2l4}-C_{3l1}+C_{4l2})].
\end{array}
$$
Under this notation, we have
$$
\begin{array}{c}
Tr_k(X^h_{\sigma})=[1+(-1)^k]\displaystyle{\frac{s(p)}{144}}X(s)
+\displaystyle{\frac{1}{12}}(\frac{ts(p)}{6}-2)X(s)\\[6pt]
+\mathit{Trace}_{h_t}\,\{{\cal V}_{\sigma}\ni V\to
[\displaystyle{\frac{t}{8}}g((\nabla_{X}{\cal B})(V),{\cal B}(V))\\[6pt]
\hspace{6.5cm}+(-1)^{k+1}\displaystyle{\frac{ts(p)}{24}}g(\delta{\cal
B}(K_VX),V)]\}.\\[6pt]
+(-1)^{k+1}\displaystyle{\frac{1}{2}}\Sigma(E_1,...,E_4).
\end{array}
$$
In particular, the sum $\Sigma(E_1,...,E_4)$ does not depend on
the choice of the oriented orthonormal basis $E_1,...,E_4$
(clearly it does not depend on the choice of the $h_t$-orthonormal
basis $V_1,V_2$ of ${\cal V}_{\sigma}$ as well). Since
$$
\Sigma(E_3,E_4,E_1,E_2)=-\Sigma(E_1,E_2,E_3,E_4),
$$
it follows that
$$
\Sigma(E_1,E_2,E_3,E_4)=0.
$$
This proves the lemma.


\begin{thebibliography}{1000}


\bibitem{AHS} M.F. Atiyah, N.J. Hitchin, I.M. Singer, {\it Self-duality in four-dimensional Riemannian geometry},
Proc. Roy. Soc. London, Ser.A 362 (1978), 425-461.

\bibitem{Besse} A. Besse, {\it Einstein manifolds}, Classics in Mathematics, Springer-Verlag, Berlin, 2008.

\bibitem{BLS} G.~Bor, L.~Hern\'andez-Lamoneda, M.~Salvai, {\it
Orthogonal almost-complex structures of minimal energy}, Geom.
Dedicata {\bf 127} (2007), 75-85.

\bibitem{CG} E.~Calabi, H.~Gluck, {\it What are the best almost-complex structures on the 6-sphere?},
Proc. Sym. Pure Math. {\bf 54} (1993), part 2, 99-106.


\bibitem{D05} J. Davidov, {\it Einstein condition and twistor spaces of compatible partially
complex structures}, Diff. Geom. and its Appl. {\bf 22} (2005),
159-179

\bibitem{D16}J. Davidov, {\it Harmonic almost Hermitian structures}, arXiv:1605.06804v3 [math.DG] 13 Jun 2016.


\bibitem{DM91} J. Davidov, O. Mushkarov, {\it On the Riemannian curvature of a twistor space},  Acta Math. Hungarica
{\bf 58} (1991), 319-332.

\bibitem{DM02} J. Davidov, O. Mushkarov, {\it Harmonic almost-complex structures on twistor spaces}, Israel J. Math.
{\bf 131} (2002), 319-332.

\bibitem{DHM15} J.~Davidov, A.~Ul~Haq, O.~Mushkarov, {\it Almost complex structures that are harmonic maps},
arXiv:1504.01610v2 [math.DG] 19 Aug 2015.

\bibitem{DGM} J. Davidov, G. Grantcharov, O. Mushkarov, {\it Twistorial
examples of $\ast$-Einstein manifolds}, Ann. Glob. Anal. Geom. {\bf
20} (2001), 103-115.

\bibitem{Der} A. Derdzinski, {\it Examples de metrques de K\"ahler et d'Einstein autoduales sur le plan complexe},
in: L. Berard-Bergery, M. Berger and C. Houzel (eds), Geometrie
riemannienne en dimension 4, Seminaire Arthur Besse, CEDIC/Fernand
Nathan, Paris, 1981, 334–346.

\bibitem{EL} J.~Eells, L.~Lemaire, {\it Selected topics in harmonic maps}, Cbms Regional Conference
Series in Mathematics, vol. {\bf 50}, AMS, Providernce, Rhode
Island, 1983.


\bibitem{ES} J. Eells, S. Salamon, {\it Twistorial constructions of harmonic maps of surfaces into four-manifolds},
Ann. Scuola Norm. Sup. Pisa, ser.IV, {\bf 12} (1985), 589-640.

\bibitem{FK}Th. Friedrich, H. Kurke, {\it Compact four-dimensional self-dual Einstein
manifolds with positive scalar curvature}, Math. Nachr. {\bf 106}
(1982), 271-299.

\bibitem{G} A. Gray, {\it Minimal varieties and almost Hermitian
manifolds}, Michigan Math. J. {\bf 12} (1965), 273-287.

\bibitem{H81} N.J. Hitchin, {\it K\"ahlerian twistor spaces}, Proc. London Math. Soc. III
Ser. {\bf 43} (1981), 133-150.

\bibitem{H95} N. Hitchin, {\it Twistor spaces, Einstein metrics and isomonodromic deformations}, J. Differential
Geom. {\bf 42} (1995), 30–112.

\bibitem{K} T. Kato, {\it Perturbation theory for linear operators},
Springer-Verlag, Berlin-Heidelberg-New York, 1980.


\bibitem{L} P. Lax, {\it Linear algebra and its applications},
John \& Sons, Inc., Hoboken, New Jersey, 2007.

\bibitem{Le1} C. LeBrun, {\it H-space with a cosmological constant}, Proc. Roy. Soc.
London A {\bf 380} (1982), 171–185.

\bibitem{Le2} C. LeBrun, {\it Counter-example to the generalized positive action
conjecture}, Comm. Math. Phys. {\bf 118} (1988), 591–596.

\bibitem{Le3} C. LeBrun, {\it Explicit self-dual metrics on ${\mathbb C}{\mathbb P}^2\sharp . . .
\sharp {\mathbb C}{\mathbb P}^2$}, J. Differential Geom. {\bf 34}
(1991), 223– 253.

\bibitem{M} O. Mu\v skarov, {\it  Structures presque hermitienes sur espaces
twistoriels et leur types}, C. R. Acad. Sci. Paris S\'er.I Math.
{\bf 305} (1987), 307-309.

\bibitem{Ped} H. Pedersen, {\it Einstein metrics, spinning top
motions and monopoles}, Math. Ann. {\bf 274} (1986), 35–59.

\bibitem{R} F. Rellich, {\it Perturbation theory of eigenvalue
problems}, Notes on mathematics and its applications, Gordon and
Breach science publishers, New York-London-Paris, 1969.

\bibitem {Tod} K. P. Tod, {\it The $SU(\infty)$-Toda field equation and special
four-dimensional metrics}, In: J. E. Andersen, J. Dupont, H.
Pedersen and A. Swann (eds), Geometry and Physics, Aarhus, 1995,
Lecture Notes in Pure Appl. Math. {\bf 184}, Marcel Dekker, New
York, 1997, 307-312.


\bibitem{V} J. Vilms, Totally geodesic maps, J. Diff. Geom. 4 (1970), 73-79.

\bibitem{W1}C. M. Wood, {\it Instability of the nearly-K\"ahler six-sphere}, J. reine angew. Math.
{\bf 439} (1993), 205-212.


\bibitem{W2} C. M. Wood, {\it Harmonic almost-complex structures}, Compositio Mathematica {\bf 99}
(1995), 183-212.



\end{thebibliography}
\end{document}